\documentclass{amsart}
\usepackage{amssymb, amsfonts}
\usepackage{xypic}
\usepackage{amsthm}

\newtheorem{thm}{Theorem}[section]
\newtheorem{cor}[thm]{Corollary}
\newtheorem{pro}[thm]{Proposition}
\newtheorem{lem}[thm]{Lemma}

\theoremstyle{definition}
\newtheorem{defn}[thm]{Definition}

\newtheorem{exmp}[thm]{Example}

\newtheorem{notn}[thm]{Notation}

\newtheorem{rem}[thm]{Remark}

\newcommand{\p}{\text{pro-}}

\newcommand{\mm}{\mathcal{M}}

\newcommand{\ee}{\mathcal{P}}
\newcommand{\pp}{\mathcal{P}}

\newcommand{\hh}{\mathcal{H}}

\newcommand{\dd}{\mathcal{D}}
\newcommand{\cc}{\mathcal{C}}

\newcommand{\rarr}{\rightarrow}
\newcommand{\map}{\rightarrow}

\newcommand{\holim}{\text{holim}}

\newcommand{\colim}{\text{colim}}

\newcommand{\zz}{\mathbb{Z}}
\newcommand{\nn}{\mathbb{N}}
\newcommand{\ess}{\text{essentially levelwise}}

\newcommand{\Ho}{\text{Ho}}
\newcommand{\Map}{\text{Map}}
\renewcommand{\lim}{\text{lim}}

\newcommand{\fib}{\text{hofib}\,}
\newcommand{\hofib}{\text{hofib}\,}

\newcommand{\hocofib}{\text{hocofib}\,}
\newcommand{\co}{\text{co}}
\newcommand{\im}{\text{im}}
\newcommand{\inj}{\text{inj-}}
\newcommand{\proj}{\text{proj-}}
\newcommand{\col}{\colon\thinspace}
\makeatletter
\let\c@equation\c@thm
\makeatother \numberwithin{equation}{section}

\newcommand{\dfn}{\textbf} %Make defined words bold.
\newcommand{\mdfn}[1]{\dfn{\mathversion{bold}#1}}
\title{T-model structures}

\author{Halvard Fausk}
\author{Daniel C.~Isaksen}

\address{Department of Mathematics, University of Oslo,
1053 Blindern, 0316 Oslo, Norway}
\address{Department of Mathematics, Wayne State University,
Detroit, MI 48202, USA}

\email{fausk@math.uio.no} \email{isaksen@math.wayne.edu}

\date{October 9, 2006}

\begin{document}
%\tableofcontents

\begin{abstract}
For every stable model category $\mm$ with a certain extra
structure, we produce an associated model structure on the
pro-category $\p \mm$ and a spectral sequence, analogous to the
Atiyah-Hirzebruch spectral sequence, with reasonably good
convergence properties for computing in the homotopy category of $\p
\mm$. Our motivating example is the category of pro-spectra.

The extra structure referred to above is a t-model structure. This
is a rigidification of the usual notion of a t-structure on a
triangulated category. A t-model structure is a proper simplicial
stable model category $\mm$ with a t-structure on its homotopy
category together with an additional factorization axiom.
\end{abstract}

\maketitle

\section{Introduction}

Recent efforts to understand the homotopy theory of pro-objects have
resulted in several different model structures on pro-categories,
such as the strict model structure \cite{EH} \cite{strict}, the $\pi_*$-model
structure on pro-spaces \cite{pro}, and the $\pi_*$-model structure
on pro-spectra \cite{ipi}. The arguments required for establishing
these model structures are similar, yet the published proofs are
distinct and have an ad hoc flavor. In the accompanying paper
\cite{ffi} we develop a general framework of filtered model
categories for giving model structures on pro-categories.
%This
%highly technical theory of filtered model categories
%ystematize the constructions of model structures on
%pro-categories and contains all examples we are aware of.

In this paper we explore a particular class of filtered model
structures on {\em stable} model categories. These filtered model
structures arise from t-structures on the homotopy category of the
stable model category (recall that such a homotopy category is a
triangulated category).

More precisely, we work with a {\em t-model structure}. This is a
proper simplicial stable model category with a t-structure on its
homotopy category together with a certain kind of lift of the
t-structure to the model category itself. This ``rigidification'' of
the t-structure is expressed in terms of an additional factorization
axiom for the model category.
%Our experience with examples is too limited to we expect
%that many common t-structures can be rigidified in this
%fashion.
%See Proposition \ref{cofibrantlygenerated} for a
%concrete result along these lines.

If $\mm$ is a t-model category, then we produce a model structure on
the category $\p \mm$ (see Theorem \ref{thm:pro-ms}). We show that
the associated homotopy category $\ee$ of $\p \mm$ has a t-structure
(see Proposition \ref{pro:tprostructure}). One important property of
this t-structure is that an object lies in $\pp_{\geq n}$ for all
$n$ if and only if it is contractible (see Proposition
\ref{pro:tprostructure} again). This property has at least two
important consequences. First, it allows us to construct an
Atiyah-Hirzebruch type spectral sequence with reasonable convergence
properties (see Theorem \ref{thm:pro-specseq}). Second, it allows us
to prove results for detecting weak equivalences analogous to the
Whitehead theorem (see Theorem \ref{thm:pro-hh}).

%The following result
% summarizes the main
%properties of the model structure on the pro-category
%associated to a t-model structure.
%
%\begin{thm} Let $\mm$ be a t-model category. There is a
%proper simplicial stable model structure on $\p \mm$. There is
%a levelwise defined t-structure $\ee^{\leq n} , \ee^{\geq n}$
%on the homotopy category $\ee$ of $\p \mm$. The weak
%equivalences on $\p \mm$ are so that an object $ X \in \ee$ is
%isomorphic to the null object if and only if $ X \in \cap_n
%\ee^{\leq n}$.
%\end{thm}
%
%The property that $ X \in \p \mm$ is isomorphic to the null
%object if and only if $ X \in \cap_n \ee^{\leq n}$ alow us to
%set up an Atiyah-Hirzebruch type spectral sequence with
%reasonable convergence properties \ref{thm:pro-specseq}, and
%to prove Whitehead type results for detecting weak
%equivalences \ref{thm:pro-hh}.

Although our main interest in t-model structures is to produce model
structures on pro-categories, the notion of a t-model structure is
likely to be useful in other contexts. Motivic homotopy theory
\cite{MV} is a combination of ideas from homotopy theory and from
motivic algebraic geometry. Since model categories are important in
homotopy theory and since t-structures are important in the study of
motives, the interaction between these two notions is probably
relevant in motivic homotopy theory. See \cite{Morel} for a possible
example.

%We consider two main examples of t-model structures. First,
%there is a t-model structure for spectra that is associated to
%the t-structure on the stable homotopy category given by
%Postnikov sections. Second,
%there is a t-model structure for chain complexes of
%$R$-modules; it is associated to the t-structure for
%chain complexes.

One specific example of a model structure on a pro-category obtained
from a t-model structure is the $\pi_*$-model structure for
pro-spectra. It is obtained from any reasonable model category of
spectra, where the t-structure on the homotopy category of spectra
is given by Postnikov sections. The original motivation for this
paper was an extension of this model structure to a category of
pro-$G$-spectra when $G$ is a profinite group. That extension is
treated in a separate paper \cite{sfi}.

Another example is the $H_*$-model structure on the category of
pro-chain complexes of modules over a unital ring, in which weak
equivalences are detected by pro-homology groups.
This model structure for pro-chain complexes is obtained
from the projective model structure on the category of chain
complexes equipped with the standard t-structure on its derived
category.
%The model structure we produce on the category of pro-chain complexes
%has a homotopy category which is itself triangulated and has a well behaved
%t-structure. Hence
%is a more natural object to study then the pro-derived
%category which does not have a triangulated structure.
%In both of these examples, the t-structure on the homotopy category
%lifts to a t-model structure.

%We consider only these two simple minded examples of t-model
%structures
% in this
%paper as they already give interesting model structure on the
%pro-category. A summary of our main results applied to the
%example of spectra with a t-structure given by Postnikov
%sections is given below. We do not make the results applied to
%the case of the pro-chain complexes explicit.
%The details of the corresponding example of chain
%complexes is left to the reader to extract from our general
%results.

\subsection{$\pi_*$-model structure on pro-spectra}

We provide a summary of our results for the $\pi_*$-model structure
on pro-spectra. There are similar results for the $H_*$-model
structure on pro-chain complexes.
% cohomology functor in Theorem \ref{thm:pro-specseq} is then a suitable
%generalization of the hyperext functor.

We remind the reader that the theorems that appear later are much
more general. The pro-categorical terminology is established in
Section \ref{sctn:pro-prelim}.

\begin{defn}
A map $f$ of pro-spectra is a $\pi_*$-weak equivalence if:
\begin{enumerate}
\item
$f$ is an essentially levelwise $m$-equivalence for some $m$, and
\item
$\pi_n f$ is a pro-isomorphism of pro-groups for all integers $n$.
\end{enumerate}
\end{defn}

We acknowledge that the first condition above appears unnatural at
first glance. However, we suspect that it is not possible to
construct a model structure on pro-spectra if this condition is not
included. In fact, the definition that appears later is different
(see Definition \ref{defn:pro-ms}). Here we have given a more
concrete equivalent reformulation (see Theorem \ref{thm:pro-hh}).

\begin{thm}[Theorem \ref{thm:pro-ms}]
There is a model structure on the category of pro-spectra in which
the weak equivalences are the $\pi_*$-weak equivalences.
\end{thm}

The cofibrations in this model structure are, up to isomorphism,
the levelwise cofibrations. The fibrations are more complicated to
describe. Section \ref{sctn:pro-model} contains a reasonably
concrete characterization of the fibrations.

One of the key observations about this model structure is that the
map $X \map \{ P_n X \}$ from a spectrum to its Postnikov tower is a
$\pi_*$-weak equivalence. In fact, Postnikov towers are the key
ingredient in constructing fibrant replacements.

One of the main uses of the previous theorem is the construction of
an Atiyah-Hirzebruch spectral sequence for pro-spectra.

\begin{thm}[Theorem \ref{thm:pro-specseq}]
Let $X$ and $Y$ be pro-spectra. There is a spectral sequence with
\[
E_2^{p,q} = H^p(X; Y^q).
\]
The spectral sequence converges conditionally to $[X,Y]^{p+q}$ if:
\begin{enumerate}
\item
$X$ is uniformly bounded below (i.e., there exists an integer $N$
such that $\pi_n X_s = 0$ for all $n \leq N$ and all $s$), or
\item
$X$ is essentially levelwise bounded below (i.e., for each $s$,
there exists an integer $N$ such that $\pi_n X_s = 0$ for $n \leq
N$) and $Y$ is a constant pro-spectrum.
\end{enumerate}
\end{thm}

In the previous theorem, the notation $[X,Y]^{p+q}$ refers to weak
homotopy classes of maps of degree $p+q$ from $X$ to $Y$ in the
homotopy category of pro-spectra.
The $E_2$-term $H^p(X; Y^q)$
is singular cohomology of the pro-spectrum $X$ with coefficients
in the pro-abelian group $Y^q = \pi_{-q} Y$. Recall that the $p$-th
cohomology group $H^p(X; A)$ of a pro-spectrum $X$ with coefficients
in an abelian group $A$ is defined to be $\colim_s H^p(X_s;A)$. The
definition in the general case is obtained from Definition
\ref{Ecohomology} and Propositions \ref{cor:homset} and
\ref{pro:rigid-heart}.
%The $p$th cohomology group $H^p(X; \{ A_a \})$ of a
%pro-spectrum $X$ with coefficients in a pro-abelian group $\{
%A_a \}$ is defined to be
%\[ \pi_0 ( \holim_a \colim_s \Map ( X_s , \Sigma^p H A_a ) )
%\] where $\{ HA_a \}$ is a functorial fibrant
% Eilenberg Mac\, Lane spectrum associated to $A_a$,
% and each $X_s$ are cofibrant.
%There is cohomological Whitehead theorem \ref{thm:Whitehead}.

\subsection{Summary}
We summarize the contents of the paper by section.

We give a short review of the basic properties of t-structures in
Section \ref{t}, assuming no prior knowledge of t-structures. In
Section \ref{subsctn:n-eq}, starting with a t-structure on the
homotopy category of a stable model category $\mm$, we produce a
filtration on the class of morphisms in $\mm$. We reformulate some
basic properties of t-structures in this language. In Section
\ref{sctn:stable} we introduce t-model categories.

 The rest of the paper is concerned with pro-categories.
The basic theory of pro-categories is briefly reviewed in Section
\ref{sctn:pro-prelim}. We also review the strict model structure and
discuss its mapping spaces. In Section \ref{sctn:pro-model} we show
that a t-model category gives rise to a filtered model category, and
we use this to give a model structure on its pro-category. For
reasons that will be apparent later, a model structure on a
pro-category obtained in this way is called an $\hh_*$-model
structure. We describe the cofibrations and fibrations of
$\hh_*$-model structures in some detail and discuss Quillen
equivalences between pro-categories with $\hh_*$-model structures.
%model structures on pro-categories
%obtained from t-model structures.
We then introduce functorial towers of truncation functors in
Section \ref{sctn:Postnikov}. They are used to form fibrant
replacements and also to construct the Atiyah-Hirzebruch spectral
sequence. We next describe the weak homotopy type of mapping spaces
in $\hh_*$-model structures in Section \ref{sctn:homotopy-class}.
This is used in Section \ref{sctn:t-structureonprocategoreis} to
give a t-structure on the homotopy category of an $\hh_*$-model
structure.
%on a pro-category associated to t-model structure.
Under reasonable assumptions, we identify the heart of this
t-structure. In Section \ref{ah} we construct an Atiyah-Hirzebruch
spectral sequence for triangulated categories with a t-structure.
 The
spectral sequence has reasonably good convergence properties when
applied to the t-structure on the homotopy category of an
 $\hh_*$-model structure.

The last two sections of the paper are devoted to multiplicative
structures on pro-categories.
%%%One reason for including this material is to
%We warn the reader that tensor structures on pro-categories do
%not have all the properties one might hope for.
In Section
\ref{sec:protensor} we give some basic results concerning the
interaction of tensor structures and pro-categories. In Section
\ref{sec:protensormodelstructure} we discuss tensor model categories
% The $\hh_*$-model structures on
%pro-categories do not behave well with respect to tensor structures.
%In general we can not
%expect to get a well behaved tensor product on the homotopy category.
and show that we get a partially defined tensor product for some
objects in the homotopy category.
%The strict model structure is far
%better behaved; in general there is a induces tensor structure
%on its homotopy category.
At the very end, we consider multiplicative structures on the
Atiyah-Hirzebruch spectral sequence constructed in Section \ref{ah}.

\subsection{Conventions}

We assume that the reader is familiar with model categories.
The reference \cite{hov} is particularly relevant because of its emphasis
on stable model categories, and \cite{hir} is also suitable.

We also assume that the reader has a certain practical familiarity
with pro-categories. Although a brief review is given in Section
\ref{sctn:pro-prelim}, see \cite{AM}, \cite{EH}, or \cite{pro} for
additional background.

We use homological grading when working with triangulated categories
and t-structures. This disagrees with the more common cohomological
grading (see \cite{bbd} for example), but it is more convenient
from the perspective of stable homotopy theory.
To emphasize the notational distinction, we use lower subscripts instead
of upper subscripts.

%%%We have decided to break the usual convention in homological algebra
%%%which is to work with cohomological grading \cite{Manin}.
%%%We warn the reader that this is somewhat unnatural from the
%%%perspective of stable homotopy theory, where the usual grading is
%%%homological.
%%%Thus, a translation is required when applying the general results of
%%%this paper to the specific example of spectra. For example, the
%%%usual notion of an ``$n$-equivalence'' of spectra corresponds to our
%%%notion of a cohomological $(-n)$-equivalence. Also, the usual notion
%%%of a ``bounded above'' spectrum corresponds to our notation of a
%%%bounded below object.

Throughout the paper, $\mm$ is a proper simplicial stable model
category. We always assume that $\mm$ has functorial factorizations,
even though the model structure on $\p \mm$ does {\em not} necessarily have
functorial factorizations. We let $\dd$ stand for the homotopy
category of $\mm$ because the notation $\Ho(\mm)$ is too cumbersome
for our purposes. Finally, $\ee$ stands for the homotopy category of
$\p \mm$.

The simplicial assumption on $\mm$ is probably not necessary for any
of our main results, but it is a very convenient hypothesis. Most of
the results go through with a weakening of this assumption; see
\cite[Sec.~6]{fch} for more details.
On the other hand, the properness
assumption on $\mm$ is essential for the existence of
model structures on $\p \mm$; see \cite{strict} for an explanation.

If $\cc$ is any category containing objects $X$ and $Y$, then
$\cc(X,Y)$ denotes the set of morphisms from $X$ to $Y$.
Occasionally we will use the notation $[X,Y]$ for the set of
morphisms in a homotopy category; in this case, the context makes
the precise meaning clear.

\section{t-structures}
\label{t}

In this section we give a short review of the theory of t-structures
on triangulated categories. We only discuss the aspects that are
relevant for our work. The original source for this material is
\cite{bbd}, but we refer to \cite{Manin} whenever possible. We
assume the reader has a working knowledge of triangulated
categories.

Let $\dd$ be a triangulated category, and let $\Sigma$ denote the
shift functor so that distinguished triangles are of the form $ A
\rarr B \rarr C \rarr \Sigma A$.

\begin{defn}
\label{defn:t} A \bf t-structure \rm on $\dd$ consists of two
strictly full subcategories $\dd_{\geq 0} $ and
 $\dd_{\leq 0} $ such that
 \begin{enumerate}
 \item $\dd_{\geq 0} $ is closed under $\Sigma$, and
$\dd_{\leq 0} $ is closed under $\Sigma^{-1}$.
\item For every object $X$ in $\dd$, there is a distinguished
triangle \[ X' \rarr X \rarr X'' \rarr \Sigma X' \] such that $X'
\in \dd_{\geq 0} $ and $X'' \in \Sigma^{-1} \dd_{\leq 0}$.
\item $\dd ( X ,Y ) = 0$
whenever $ X \in \dd_{\geq 0} $ and $Y \in \Sigma^{-1} \dd_{\leq 0}$.
\end{enumerate}
\end{defn}

The reader who is already familiar with t-structures should keep in mind
that we are using homological grading, not cohomological grading.

Recall that a subcategory is strictly full if it is full and if it
is closed under isomorphisms.

In any t-structure, there are two $\zz$-graded families of strictly
full subcategories defined by
$\dd_{\geq n} = \Sigma^{n} \dd_{\geq 0}$ and
$\dd_{\leq n} = \Sigma^{n} \dd_{\leq 0} $ for any integer
$n$. By part (1) of Definition \ref{defn:t}, the categories
$\dd_{\geq n}$ are a decreasing filtration in the sense that
$\dd_{\geq n+1}$ is contained in $\dd_{\geq n}$, and the categories
$\dd_{\leq n}$ are an increasing filtration in the sense that
$\dd_{\leq n}$ is contained in $\dd_{\leq n+1}$. Note that parts (2)
and (3) of Definition \ref{defn:t} can now be rewritten in terms of
$\dd_{\leq -1}$ instead of $\Sigma^{-1} \dd_{\leq 0}$.

\begin{exmp}
\label{ex:standard} Let $\dd$ be the derived category of chain
complexes of modules over a ring. The \mdfn{standard t-structure} on
$\dd$ \cite[IV.4.3]{Manin} is given by
\[ \dd_{\geq n} = \{ X \ | \ H_i ( X) =0 \text{ for } i < n \}
\]
\[ \dd_{\leq n} = \{ X \ | \ H_i ( X) =0 \text{ for } i > n \} .
\]
The shift functor is defined by $(\Sigma X)_n =
X_{n-1}$, and the differential $(\Sigma X)_{n+1} \map (\Sigma X)_{n}$
is equal to the negative of the differential $X_n \map X_{n-1}$.
\end{exmp}

\begin{exmp}
\label{ex:Postnikov} Let $\dd$ be the homotopy category of spectra.
The \mdfn{Postnikov t-structure} on $\dd$ is
 given by
\[ \dd_{\geq n} = \{ X \, | \, \pi_i (X ) = 0 \text{ for } i < n \}
\]
\[
\dd_{\leq n} = \{ X \, | \, \pi_i (X ) = 0 \text{ for } i > n \}.
\]
The proof of this fact is classical stable homotopy theory. The main
point is to show that $[X,Y] = 0$ when $X$ is $(-1)$-connected
(i.e., $\pi_i X = 0$ for $i < 0$) and $Y$ is $0$-coconnected (i.e.,
$\pi_i Y = 0$ for $i \geq 0$). See
\cite[Prop.~3.6]{Margolis}, for example.
\end{exmp}
%%%At first glance, Examples \ref{ex:standard} and \ref{ex:Postnikov}
%%%may appear rather different. However, the main difference is that
%%%the second example is homologically graded rather than
%%%cohomologically graded. In other words, the $n$th stable homotopy
%%%group functor $\pi_n$ for spectra corresponds to the $(-n)$th
%%%cohomology functor $H^{-n}$ for chain complexes.
%%%

The two strictly full subcategories $ \dd_{\geq n}$ and
$\dd_{\leq n-1}$ of a t-structure determine each other as follows.

\begin{lem} \label{lem:determine}
 An object $X $ is in $ \dd_{\geq n}$ if and only
if $ \dd ( X , Y) =0$ for all $Y$ in $\dd_{\leq n - 1}$, and an
object $Y$ is in $\dd_{\leq n- 1} $ if and only if
 $ \dd ( X , Y) =0$ for all $X$ in $\dd_{\geq n}$.
\end{lem}

\begin{proof}
We prove the first claim; the proof for the second claim is similar.
One direction follows immediately from part (3) of Definition
\ref{defn:t}.

For the other direction, suppose that $\dd(X,Y) = 0$ for all $Y$ in
$\dd_{\leq n-1}$. Part (2) of Definition \ref{defn:t} says that we
can find a distinguished triangle $X' \map X \map X'' \map \Sigma
X'$ such that $X'$ belongs to $\dd_{\geq n}$ and $X''$ belongs to
$\dd_{\leq n-1}$. Now apply $\dd ( - , X'')$ to obtain a long exact
sequence. Since $X'$ and $\Sigma X'$ are both in $\dd_{\geq n}$, it
follows that $\dd(X, X'')$ and $\dd(X'', X'')$ are isomorphic. But
our assumption implies that the first group is zero, so the second
group is also zero. Thus $X''$ is isomorphic to $0$, and $X' \rarr
X$ is an isomorphism.
\end{proof}

The next corollary follows immediately from Lemma
\ref{lem:determine}.

\begin{cor}
\label{cor:t-prop} \mbox{}
\begin{enumerate}
\item
The zero object $0$ belongs to both $\dd_{\geq n}$ and $\dd_{\leq n}$
for every $n$.
\item
An object that is both in $\dd_{\geq n} $ and in $\dd_{\leq n-1}$ is
isomorphic to $0$.
\item
The subcategories $\dd_{\geq n}$ and $\dd_{\leq n }$ of $\dd$ are
closed under retract.
\end{enumerate}
\end{cor}

\begin{cor} \label{t2outof3}
Let $X \rarr Y \rarr Z \rarr \Sigma X$ be a distinguished triangle.
\begin{enumerate}
\item If $X$ and $Z$ belong to $\dd_{\geq n}$, then so does $Y$.
\item If $X$ and $Z$ belong to $\dd_{\leq n}$, then so does $Y$.
\item If $X$ belongs to $\dd_{\geq n-1}$ and $Y$ belongs to $\dd_{\geq n}$,
then $Z$ belongs to $\dd_{\geq n}$.
\item If $X$ belongs to $\dd_{\leq n-1}$ and $Y$ belongs to $\dd_{\leq n}$,
then $Z$ belongs to $\dd_{\leq n}$.
\item If $Y$ belongs to $\dd_{\geq n}$ and $Z$ belongs to $\dd_{\geq n+1}$,
then $X$ belongs to $\dd_{\geq n}$.
\item If $Y$ belongs to $\dd_{\leq n}$ and $Z$ belongs to $\dd_{\leq n+1}$,
then $X$ belongs to $\dd_{\leq n}$.
\end{enumerate}
\end{cor}

\begin{proof}
The first two claims follow immediately from Lemma
\ref{lem:determine}. The other claims follow from the first two; we
illustrate with the fifth statement.

We have an exact triangle
\[
\Sigma^{-1} Z \rarr X \rarr Y \rarr Z
\]
in which both $\Sigma^{-1} Z $ and $Y$ belong to $\dd_{\geq n}$. The
result now follows from the first claim.
\end{proof}

We now give a key lemma about t-structures.

\begin{lem}
\label{truncation} Let $X' \rarr X \rarr X'' \rarr \Sigma X' $ and $
Y' \rarr Y \rarr Y'' \rarr \Sigma Y'$ be two distinguished triangles
with $X'$ and $Y'$ in $\dd_{\geq n }$ and $X''$ and $Y''$ in
$\dd_{\leq n -1}$. Let $f \col X \rarr Y $ be a map in $\dd$. Then
there are unique maps $ f' \col X' \rarr Y'$ and $ f'' :X'' \rarr
Y''$ such that there is a commutative diagram
\[ \xymatrix{
X' \ar[r] \ar[d]_{f'} & X \ar[d]_f \ar[r] & X'' \ar[d]_{f''} \ar[r]
&
 \Sigma X' \ar[d]^{\Sigma f'} \\
Y' \ar[r] & Y \ar[r] & Y'' \ar[r] & \Sigma Y'. }
\]
\end{lem}

\begin{proof}
By part (3) of Definition \ref{defn:t}, we have that $\dd ( X' , Y''
) = 0$. Hence there exists a map $ f' \col X' \rarr Y'$ making the
left square commute, and then also a map $ f'':X'' \rarr Y''$ such
that we get a map of distinguished triangles.
%This shows that $f'$ and $f''$ exist.

Now we have to show that $f'$ and $f''$ are unique. Let $j$ be the
map $Y' \map Y$. There is an exact sequence
\[
\xymatrix@1{ \dd ( X' , \Sigma^{-1} Y'' ) \ar[r] & \dd ( X' , Y' )
\ar[r]^{j_*} & \dd ( X' , Y ). }
\]
The left group is trivial by part (3) of Definition \ref{defn:t}
since $ \Sigma^{-1} Y''$ belongs to $\dd_{\leq n-1 }$. Therefore,
$j_*$ is injective, so $f'$ is unique. A similar argument involving
the map $X \map X''$ shows that $f''$ is also unique.
\end{proof}

The importance of Lemma \ref{truncation} is expressed in the
following proposition.

\begin{pro}
\label{pro:truncation} There are
\mdfn{truncation functors $\tau_{\geq n}$}
and \mdfn{$\tau_{\leq n}$} from $\dd$ into $\dd_{\geq n}$ and
$\dd_{\leq n}$ respectively
together with natural
transformations $\epsilon_n \col \tau_{\geq n} \rarr 1 $, $ \eta_n
\col 1 \rarr \tau_{\leq n } $, and $\tau_{\leq n -1 } \rarr \Sigma
\tau_{\geq n}$ such that
\[
\tau_{\geq n} X \rarr X \rarr \tau_{\leq n -1} X \rarr \Sigma
\tau_{\geq n} X
\]
is a distinguished triangle for all $X$. Up to canonical
isomorphism, these properties determine the truncation functors
uniquely.
\end{pro}

\begin{notn}
We usually write \mdfn{$X_{\geq n}$} and \mdfn{$X_{\leq n}$} for
$\tau_{\geq n}X$ and $\tau_{\leq n}X$ respectively.
\end{notn}

The functors $\tau_{\geq n}$ and $\tau_{\leq n}$ enjoy many useful
properties. Most of the claims in the next two paragraphs are proved
in \cite[Sec.~IV.4]{Manin}; the rest follow easily. In any case,
they are easily verifiable directly for Examples \ref{ex:standard}
and \ref{ex:Postnikov}.

The functors $\tau_{\geq n}$ and $\tau_{\leq m}$ commute (up to
natural isomorphism) for all $n$ and $m$. If $m \geq n$, then
$\tau_{\geq n} \tau_{\geq m}$ and $\tau_{\geq m} \tau_{\geq n}$ are
both naturally isomorphic to $\tau_{\geq m}$, while $\tau_{\leq n}
\tau_{\leq m}$ and $\tau_{\leq m} \tau_{\leq n}$ are both naturally
isomorphic to $\tau_{\leq n}$. Also, $\Sigma^n \tau_{\geq 0}$ is
naturally isomorphic to $\tau_{\geq n} \Sigma^n$, and
$\Sigma^n \tau_{\leq 0}$ is naturally isomorphic to
$\tau_{\leq n} \Sigma^n$.
Both $\tau_{\leq n-1 } \tau_{\geq n} $ and $\tau_{\geq n }
\tau_{\leq n-1 } $ are isomorphic to the zero functor.

\subsection{Hearts and cohomology} \label{subsctn:heart}

\begin{defn}
\label{defn:heart} The \mdfn{heart $\hh(\dd)$} of a t-structure
$\dd$ is the full subcategory $\dd_{\geq 0 } \cap \dd_{ \leq 0 }$ of
$\dd$.
\end{defn}

%In general, the intersection of full subcategories $\cc_a$ of
%a category $\cc$ is the full subcategory whose objects are
%simultaneously isomorphic to object in $\cc_a$ for all a.

For any t-structure, the heart $\hh(\dd)$ is an abelian category
\cite[Sec.~IV.4]{Manin}.

\begin{defn} \label{homology}
The \mdfn{$n$-th homology functor $\hh_n$} associated to a
t-struct\-ure is defined to be the functor $ \tau_{\leq 0}
\tau_{\geq 0} \Sigma^{-n} $.
\end{defn}

The homology functor $\hh_n$ is a covariant functor $\dd \map \hh
(\dd)$. The following lemma is proved in
\cite[Thm.~IV.4.11a]{Manin}.

\begin{lem}
\label{lem:exact-cohlgy}
Let $X \map Y \map Z \map \Sigma X$ be a
distinguished triangle. There is a long exact sequence
\[
\cdots \map \hh_k X \map \hh_k Y \map \hh_k Z \map \hh_{k-1} X \map
\cdots
\]
in the abelian category $\hh(\dd)$.
\end{lem}

\begin{defn}
\label{Ecohomology} Let $E$ be an object in the heart $\hh (\dd)$.
The \mdfn{$n$-th cohomology functor $H^n ( -; E)$ with
$E$-coefficients } is $ \dd (- , \Sigma^{n} E) $.
\end{defn}

The functor $H^n (- ;E)$ is a contravariant functor from $\dd$ to
the category of abelian groups. The following lemma follows
immediately from formal properties of triangulated categories.

\begin{lem}
\label{lem:Ecohomologyexact}
Let $X \map Y \map Z \map \Sigma X$ be a
distinguished triangle, and let $E$ belong to $\hh(\dd)$.
There is a long exact sequence
\[
\cdots \map H^k(X; E) \map H^k(Y;E) \map H^k(Z;E) \map H^{k+1}(X;E) \map
\cdots
\]
in the abelian category $\hh(\dd)$.
\end{lem}

\begin{exmp} Let $\dd$ be the triangulated category of chain
complexes with the standard t-structure (see Example
\ref{ex:standard}). The heart of $\dd$ is isomorphic to the category
of abelian groups. The $n$-th homology $\hh_n X$ of a chain
complex $X$ is the usual $n$-th homology
\[
\ker (X_n \rarr X_{n-1})/ \im (X_{n+1} \rarr X_n )
\]
of $X$. For any abelian group $E$, the $n$-th cohomology with
$E$-coefficients of a chain complex $X$ is the $n$-th hyperext group
$\text{Ext}^n ( X , E )$.
\end{exmp}

\begin{exmp} \label{exmp:spectraheart}
Let $\dd$ be the category of spectra with the Postnikov t-structure.
The heart of this t-structure is the full subcategory of
Eilenberg-Mac\,Lane spectra (in degree 0). This category is
equivalent to the category of abelian groups. The $n$-th homology
functor $\hh_n$ is the usual $n$th stable homotopy group functor
$\pi_{n}$. When $E$ is an abelian group, $n$-th
cohomology with $E$-coefficients is $n$-th singular cohomology with
coefficients in $E$.
\end{exmp}

%\begin{exmp}
% Let $\dd$ be the stable $G$-equivariant
%homotopy category for a finite group $G$, with the Postnikov
%% Then the Cohomology of the type $\dd ( - , H)^k$ for an object $H$ in
% the heart is exactly the $k$-th Bredon cohomology with coefficients in
% $M$. The heart is the category of coefficient systems, and the
 %heart cohomology is the ordinary homotopy coefficient system, with
 %a cohomological grading. However our singular cohomology with
%coefficients in a pro-abelian group is an abelian group not a
%pro-abelian group.
%\end{exmp}

The following lemma says that
the layers in the towers obtained from the two sequences of truncation functors
are easily described in terms of the homology functors.

\begin{lem}
\label{lem:tower-fibers} For every integer $n$ and every $X$ in
$\dd$, there are distinguished triangles
\[
X_{\geq n+1} \map X_{\geq n} \map \Sigma^{n} \hh_n X \map
\Sigma X_{\geq n+1}
\]
and
\[
\Sigma^{n} \hh_n X \map X_{\leq n} \map X_{\leq n-1} \map
\Sigma^{n+1} \hh_n X.
\]
\end{lem}

\begin{proof}
We construct the first distinguished triangle. The construction of
the second one is similar.

Start with the object $X_{\geq n}$ of $\dd$. There is a
distinguished triangle
\[
\tau_{\geq n+1} \tau_{\geq n} X \map
\tau_{\geq n} X \map
\tau_{\leq n} \tau_{\geq n} X \map
\Sigma \tau_{\geq n+1} \tau_{\geq n} X.
\]
This is equal to the distinguished triangle
\[
\tau_{\geq n+1} X \map \tau_{\geq n} X \map \Sigma^{n} \hh_n X \map
\Sigma \tau_{\geq n+1} X.
\]
\end{proof}

The following lemma tells us when homology and cohomology theories
detect trivial objects. It is also proved in \cite[IV.4.11(b)]{Manin}.

\begin{lem} \label{lem:Whitehead}
Let $\dd$ be a triangulated category equipped with a t-structure
such that $\cap_n \dd_{\geq n}$ consists only of the objects
isomorphic to $0$. Assume that $X$ is in $\dd_{\geq m}$ for some
$m$. Then the following are equivalent:
\begin{enumerate}
\item $X$ is isomorphic to $0$.
\item $\hh_n (X) =0$ for all $n$.
\item $H^n ( X ; E ) =0$ for all $n$ and
all $E $ in $\hh (\dd)$.
\end{enumerate}
\end{lem}

\begin{proof}
Condition (1) implies conditions (2) and (3), so we need to show
that either condition (2) or (3) implies condition (1).

Assuming either condition (2)
or (3), we prove by induction on $n$ that $ X_{\geq n} \rarr X$ is
an isomorphism for all $n $. Our assumption on
$\cap_n \dd_{\geq n}$ then implies condition (1).
Since $X$ is in $\dd_{\geq m}$, the natural map $X_{\geq n} \rarr X$
is an isomorphism whenever $n \leq m$; this is the base case
of the induction.

Now suppose for induction that the map $X_{\geq n } \rarr X$ is an isomorphism.
Condition (2) and the first part of Lemma \ref{lem:tower-fibers}
gives that the composition $X_{\geq n +1 } \rarr X_{\geq n } \rarr
X$ is an isomorphism.

On the other hand, condition (3) and our inductive assumption
implies that $\dd(X_{\geq n}, \Sigma^{n} \hh_n X)$ is zero. Now
apply the functor $ \dd ( - , \Sigma^{n} \hh_n X)$ to the first
triangle in Lemma \ref{lem:tower-fibers} and use part 3 of
Definition \ref{defn:t} to conclude that
\[
\dd (\Sigma^{n} \hh_n (X) , \Sigma^{n} \hh_n (X) ) \] is zero.
It follows that $\hh_n(X)$ is isomorphic to $0$. As in the previous
paragraph, this implies that $ X_{\geq n+1} \rarr X$ is an isomorphism.
\end{proof}

\section{$n$-equivalences and co-$n$-equivalences}
\label{subsctn:n-eq}

From now on, we no longer consider just triangulated categories but
rather proper simplicial stable model categories $\mm$ (with
functorial factorization). We write $\dd$ for the homotopy category
of $\mm$, which is automatically a triangulated category because
$\mm$ is stable \cite[7.1]{hov}. We briefly review the main
properties of stable model categories that we need. See
\cite[Ch.~7]{hov} for more details.

Recall that a stable model category $\mm$ is pointed. In addition to
unreduced tensors $X \otimes K$ and cotensors $\Map(K,X)$, a pointed
simplicial model category also has \mdfn{reduced tensors $X \wedge
K$} and \mdfn{reduced cotensors $\Map_*(K, X)$} for any pointed
simplicial set $K$ and any object $X$ of $\mm$.

The \mdfn{suspension} of any object $X$ of $\mm$ is defined to be $X
\wedge S^1$ \cite[6.1.1]{hov}. Note that this construction is
homotopically correct only if $X$ is cofibrant. In general, one must
first take a cofibrant replacement for $X$. In a simplicial stable
model category, suspension is a left Quillen functor that induces an
automorphism on $\dd$. Its associated right Quillen functor is
reduced cotensor with $S^1$, which is also known as \mdfn{loops}
\cite[6.1.1]{hov}.

The homotopy category of a stable model category is a triangulated
category whose shift functor $\Sigma$ is the left derived functor of
suspension. We use the symbol $\Omega$ for the right derived functor
of loops. Note that $\Omega$ induces the inverse shift $\Sigma^{-1}$
on $\dd$. Homotopy cofiber sequences and homotopy fiber sequences
are the same, and they induce the distinguished triangles in the
homotopy category.
%Homotopy cofiber sequences and homotopy fiber sequences induce
%distinguished triangles
% \[ X \rarr Y \rarr C \rarr \Sigma X
%\] and \[ \Omega Y \rarr F \rarr X \rarr \Sigma \Omega Y
%\] in the homotopy category, where the map $ X \rarr \Sigma \Omega Y$
%is the composite \[ X \rarr Y \stackrel{- \epsilon
%(Z)}{\longrightarrow} \Sigma \Omega
%Y \] of $X \rarr Y$ with minus the
%unit map of the $\Sigma$-$\Omega$ adjunction \cite[7.1.11]{hov}.
%\comment{This seemed unnecessarily complicated
%and distracting from the main point.}

Since $\mm$ has functorial factorizations, there are functorial
constructions of homotopy fibers and homotopy cofibers in $\mm$. We
write \mdfn{$\hocofib f$} and \mdfn{$\hofib f$} for the functorial
homotopy cofiber and homotopy fiber of a map $f$. Using the
properness assumption, there are natural maps $\hofib f \map X$ and
$Y \map \hocofib f$ for any map $f:X \map Y$. In $\dd$, $\hocofib f$
is isomorphic to $\Sigma \hofib f$.

These constructions induce functors on the homotopy category $\dd$
of $\mm$. The functoriality of these constructions is one of the
chief advantages of working with stable model categories rather than
just triangulated categories, where homotopy cofibers and homotopy
fibers are only defined up to non-canonical isomorphism.

We now lift the full subcategories given by a t-structure on $\dd$
to full subcategories on $\mm$.

\begin{defn}
\label{defn:m-leq-n} Let $\mm$ be a proper simplicial stable model
category whose homotopy category $\dd$ is equipped with a
t-structure. Let \mdfn{$\mm_{\geq n} $} be the full subcategory of
$\mm$ consisting of those objects whose weak homotopy types belong
to $\dd_{\geq n}$, and let \mdfn{$\mm_{\leq n}$} be the full
subcategory of $\mm$ consisting of those objects in $\mm$ whose weak
homotopy types belong to $\dd_{\leq n}$.
\end{defn}

Many properties of $\dd_{\geq n}$ and $\dd_{\leq n}$ from Section
\ref{t} carry over to the classes $\mm_{\geq n}$ and $\mm_{\leq n}$.
For example, $\mm_{\geq n}$ is closed under $\Sigma$, and $\mm_{\leq
n}$ is closed under $\Omega$. Also, an object $X$ belongs to
$\mm_{\geq n}$ if and only if $\Omega^n X$ belongs to $\mm_{\geq 0}$
(and similarly for $\mm_{\leq n}$). Finally, if $X$ belongs to
$\mm_{\geq n}$ and $Y$ belongs to $\mm_{\leq n-1}$, then $\dd(X,Y)$
is zero.

To explore the interaction between a stable model structure and a
t-structure on its homotopy category, we find it more convenient to
work with subclasses of morphisms associated to the t-structure
rather then the subclasses of objects given by the t-structure.
\begin{defn}
\label{defn:n-eq} Let $\mm$ be a proper stable model category whose
homotopy category $\dd$ is equipped with a t-structure. The class of
\mdfn{$n$-equivalences} in $\mm$ is
\[
\mathbf{W_n} = \{ f \ | \ \hofib f \in \mm_{\geq n} \} = \{ f \ | \
\hocofib f \in \mm_{\geq n +1} \}.
\]
The class of \mdfn{co-$n$-equivalences} in $\mm$ is
\[
\mathbf{co W_n} = \{ f \ | \ \hocofib f \in \mm_{\leq n} \} = \{ f \
| \ \hofib f \in \mm_{\leq n-1} \} .
\]
\end{defn}

\begin{exmp}
Consider the standard t-structure on the derived category of chain
complexes from Example \ref{ex:standard}. A map is an
$n$-equivalence if and only if it induces a homology isomorphism
in degrees strictly less than $n$ and it induces a surjection in
degree $n$. Similarly, a map is a co-$n$-equivalence if and only if
it induces a homology isomorphism in degrees strictly greater
than $n$ and it induces an injection in degree $n$.
\end{exmp}

\begin{exmp}
Let $\mm$ be a model category of spectra with the Postnikov
t-structure on its homotopy category (see Example
\ref{ex:Postnikov}). A map is an $n$-equivalence if and only if it
induces isomorphisms in homotopy groups below dimension $n$ and it
induces a surjection in dimension $n$. Similarly, a map is a
  co-$n$-equivalence if and only if it
induces isomorphisms in homotopy groups above dimension $n$ and it
induces an injection in dimension $n$.
%%%This terminology unfortunately conflicts with the usual notion of
%%%$n$-equivalences and co-$n$-equivalences of spectra. This is a
%%%result of our use of cohomological grading, rather than the
%%%homological grading that is standard for spectra.
\end{exmp}

We now translate many of the results from Section \ref{t} into
properties of $n$-equivalences and co-$n$-equivalences.
%The verifications are left to the reader.

\begin{lem}
\label{lem:W^n-prop} \mbox{}
\begin{enumerate}
\item
The class $W_n$ is closed under $\Sigma$, and the class $\co W_n$ is
closed under $\Omega$.
\item
The class $W_n$ is contained in $W_{n-1}$, and the class $\co
W_{n-1}$ is contained in $\co W_n$.
\item
The classes $W_n$ and $\co W_n$ contain all weak equivalences.
\item
If a map is both an $n$-equivalence and a co-$n$-equivalence, then
it is a weak equivalence.
\item
The classes $W_n$ and $\co W_n$ are closed under retract.
\end{enumerate}
\end{lem}

\begin{proof}
The first two claims follow from the properties of $\mm_{\geq n}$
and $\mm_{\leq n}$ stated after Definition \ref{defn:m-leq-n} and
the fact that the functors $\Sigma$ and $\Omega$ commute with the
functors $\hofib$ and $\hocofib$ up to weak equivalence.

The next two claims follow immediately from parts (1) and (2) of
Corollary \ref{cor:t-prop}, together with the fact that a map is a
weak equivalence if and only if its homotopy fiber is contractible.

For the fifth claim, let $g$ be a retract of a map $f$. Then $\hofib
g$ is a retract of $\hofib f$ because homotopy fibers are
functorial. Now we just need to use part (3) of Corollary
\ref{cor:t-prop}.
\end{proof}

\begin{lem} \label{W^n-2outof3}
Let $f \col X \rarr Y$ and $g \col Y \rarr Z$ be two maps.
\begin{enumerate}
\item If $f$ and $g$ both belong to $W_n$, then so does $gf$.
\item If $f$ and $g$ both belong to $\co W_n$, then so does $gf$.
\item If $f$ belongs to $W_{ n-1}$ and $gf$ belongs to $W_{ n}$,
then $g$ belongs to $W_{ n}$.
\item If $f$ belongs to $co W_{ n - 1}$ and $gf$ belongs to $coW_{ n}$,
then $g$ belongs to $co W_n$.
\item If $g$ belongs to $W_{n+1}$ and $gf$ belongs to $W_n$,
then $f$ belongs to $W_n$.
\item If $g$ belongs to $coW_{n+1}$ and $gf$ belongs to $coW_n$,
then $f$ belongs to $coW_n$.
\end{enumerate}
\end{lem}

\begin{proof}
We have a distinguished triangle
\[
\fib f \rarr \fib gf \rarr \fib g \rarr \Sigma \fib f
 \]
in the homotopy category $\dd$ of $\mm$. Corollary \ref{t2outof3}
gives the desired results.
\end{proof}

\begin{lem} \label{proper}
The classes $W_n$ and $\co W_n$ are both closed under base changes
along fibrations and cobase changes along cofibrations.
\end{lem}

\begin{proof}
In a proper model structure, the homotopy fiber of a map is weakly
equivalent to the homotopy fiber of its base change along a
fibration. Similarly, the homotopy cofiber of a map is weakly
equivalent to the homotopy cofiber of its cobase change along a
cofibration.
\end{proof}

The next result gives a general setting in which
homology and cohomology functors detect weak equivalences in $\mm$.

\begin{thm} \label{thm:Whitehead} Assume that $\cap_n W_n$ is equal
to the class of weak equivalences. Let $f\col X \map Y$ be an
$m$-equivalence for some $m$. The following are equivalent:
\begin{enumerate}
\item $f$ is a weak equivalence.
\item
$\hh_n ( f )$ is an isomorphism in the heart $\hh (\dd)$ for all
$n$.
\item
$H^n( Y; E) \map H^n ( X; E) $ is an isomorphism for all $n$ and all
$E$ in $\hh(\dd)$.
\end{enumerate}
\end{thm}

Applied to the Postnikov t-structure on spectra from Example
\ref{ex:Postnikov}, this is the usual Whitehead theorem for
detecting weak equivalences with stable homotopy groups or with
ordinary cohomology.

\begin{proof}
Note first that $\cap_n W_n$ is equal to the class of weak
equivalences if and only if $\cap_n \dd_{\geq n}$ consists only of
contractible objects. This follows immediately from the definitions
and the fact that a map is a weak equivalence if and only if its
homotopy fiber is contractible.

Now Lemma \ref{lem:Whitehead} gives the desired result, using the
long exact sequences of Lemmas \ref{lem:exact-cohlgy} and
\ref{lem:Ecohomologyexact}.
\end{proof}

\begin{rem}
In Theorem \ref{thm:Whitehead}, an additional assumption on the
t-structure allows one to avoid the assumption that $f$ is
an $m$-equivalence. Namely, if both $ \cap_n W_n$ and $\cap_n \co
W_n$ are equal to the class of weak equivalences, then a map $f$ is
a weak equivalence in $\mm$ if and only if $\hh_n(f)$ is an
isomorphism for all $n$; the proof of \cite[IV.4.11]{Manin} can be
easily adapted to show this. Both the standard t-structure on chain
complexes from Example \ref{ex:standard} and the Postnikov
t-structure on spectra from Example \ref{ex:Postnikov} satisfy this
condition.

Our assumptions in Theorem \ref{thm:Whitehead} are dictated by the
t-structures we consider on homotopy categories of pro-categories.
In that case we have that $\cap_n W_n$ is equal to the class of weak
equivalences, while $\cap_n \co W_n$ is never equal to
the class of weak equivalences except in trivial cases.
See Lemma \ref{rem:constant} for more details.
\end{rem}

\begin{rem}
Consider the situation of a t-structure on a triangulated category
$\dd$ that is not associated to a stable model category. Even though
homotopy fibers and homotopy cofibers are not well-defined in $\dd$,
one can still define classes of $n$-equivalences and
co-$n$-equivalences in $\dd$ as in Definition \ref{defn:n-eq}. The
point is that homotopy fibers and homotopy cofibers are well-defined
up to non-canonical isomorphism, and that is good enough for the
purposes of Definition \ref{defn:n-eq}. All of the lemmas of this
section remain true except for Lemma \ref{proper}, which does not
make sense without a model structure.
\end{rem}

\section{t-model structures} \label{sctn:stable}

We continue to work in a proper simplicial stable model category
$\mm$ whose homotopy category $\dd$ has a t-structure. We need to
assume that the t-structure on $\dd$ can be rigidified in a certain
sense.

\begin{defn} \label{defn:tmodel}
A {\bf t-model structure} is a proper simplicial stable model
category $\mm$ equipped with a t-structure on its triangulated
homotopy category $\dd$ together with functorial factorizations of
maps in $\mm$ into $n$-equi\-val\-ences followed by
co-$n$-equi\-val\-ences.
\end{defn}

There is a t-model structure on the category of chain complexes that
induces the standard t-structure of Example \ref{ex:standard}
(except, possibly, for the simplicial structure). This example is
dealt with in greater detail in \cite[4.8,6]{fch}. More importantly
for our applications, all reasonable model categories of spectra
have t-model structures that induce the Postnikov t-structure on the
stable homotopy category. To factor a map into an $n$-equivalence
followed by a co-$n$-equivalence, apply the small object argument to
the set of maps consisting of all generating acyclic cofibrations
and also generating cofibrations whose cofibers are spheres of
dimension greater than $n$.

\begin{lem}
\label{lem:rigid-truncation} The truncation functors $\tau_{\geq
n}\col \dd \map \dd_{\geq n}$ and $\tau_{\leq n} \col \dd \map
\dd_{\leq n}$ can be lifted to functors $\tau_{\geq n}\col \mm \map
\mm_{\geq n}$ and $\tau_{\leq n}\col \mm \map \mm_{\leq n}$.
Similarly, the natural transformations $\epsilon_n \col \tau_{\geq
n} \rarr 1 $ and $\eta_n \col 1 \rarr \tau_{\leq n }$ can be lifted
to natural transformations on $\mm$ such that $\epsilon_n$ is a
natural co-$(n-1)$-equivalence, $\eta_n$ is a natural
$(n+1)$-equivalence, and $\tau_{\geq n} X
\stackrel{\epsilon_n}{\longrightarrow} X \stackrel{\eta_{n -
1}}{\longrightarrow} \tau_{\leq n - 1} X $ is a natural homotopy
fiber sequence in $\mm$.
\end{lem}

\begin{proof}
Given any object $X$ of $\mm$, factor the map $* \map X$
functorially into an $(n-1)$-equivalence $* \map X'$ followed by a
co-$(n-1)$-equivalence $X' \map X$. Define $\tau_{\geq n} X$ to be
$X'$, and $\epsilon_n(X)$ to be the natural map $X' \map X$. Define
$\tau_{\leq n}X $ to be the homotopy cofiber of $ \epsilon_{n+1}
(X)$ and $\eta_n (X)$ to be the map $ X \rarr \tau_{\leq n} X$. We
have that $X'$ belongs to $\mm_{\geq n}$ since $* \map X'$ is an
$(n-1)$-equivalence,
 and $\tau_{\leq n} X$ belongs to $\mm_{\leq n}$ since it the
co-$n$-equivalences are defined to have homotopy cofibers in $\mm_{\leq n}$.
\end{proof}

\begin{defn}
Let $\mm$ be a t-model structure. A map in $\mm$ is an
\mdfn{$n$-cofibration} if it is both a cofibration and an
$n$-equivalence. A map in $\mm$ is a \mdfn{co-$n$-fibration} if it
is both a fibration and a co-$n$-equivalence.
\end{defn}

\begin{lem} \label{cfretract}
The classes of $n$-cofibrations and co-$n$-fibrations are closed
under composition and retract.
\end{lem}

\begin{proof}
This follows immediately from the fact that the classes of
cofibrations, fibrations, $n$-equivalences, and co-$n$-equivalences
are all closed under composition and retract by part (5) of Lemma
\ref{lem:W^n-prop} and parts (1) and (2) of Lemma \ref{W^n-2outof3}.
\end{proof}

\begin{lem} \label{factor}
There is a functorial factorization of maps in $\mm$ into
$n$-cofibrations followed by co-$n$-fibrations.
\end{lem}

\begin{proof}
We construct a factorization explicitly. Let $f \col X \rarr Y$ be a
map in $\mm$. We have a natural diagram
\[\xymatrix{ X \ar[rr]^u \ar@{>->}[dr] & & Z \ar[rr]^v
\ar@{>->}[rd]^{\sim}
 & & Y \\
& A \ar@{>->}[dr]_{\sim} \ar@{>>}[ur]^{\sim}
 \ar[rr]^{\sim} & & B \ar@{>>}[ur] & \\
& & C \ar@{>>}[ur]_{\sim} & & } \] obtained as follows. First,
factor $f$ into an $n$-equivalence $u\col X \map Z$ followed by a
co-$n$-equivalence $v\col Z \map Y$. Next, factor $u$ into a
cofibration $X \map A$ followed by an acyclic fibration $A \map Z$,
and factor $v$ into an acyclic cofibration $Z \map B$ followed by a
fibration $B \map Y$. Now the composition $A \map Z \map B$ is a
weak equivalence, so it can be factored into an acyclic cofibration
$A \map C$ followed by an acyclic fibration $C \map B$.

The composition $X \rarr A \map C$ is a cofibration because it is a
composition of two cofibrations, and it is an $n$-equivalence by
Lemma \ref{W^n-2outof3}.

Similarly, the composition $C \map B \map Y$ is a co-$n$-fibration.
\end{proof}

We next prove that the classes of $n$-cofibrations and
co-$n$-fibrations determine each other via lifting properties.

\begin{lem} \label{lift}
A map is an $n$-cofibration if and only if it has the left lifting
property with respect to all co-$n$-fibrations. A map is a
co-$n$-fibration if and only if it has the right lifting property
with respect to all $n$-cofibrations.
\end{lem}

\begin{proof}
Let $i$ be an $n$-cofibration and $p$ a co-$n$-fibration. We use the
abstract obstruction theory of \cite{obs} to show that there exists
a lift $B \rarr X$ in the diagram
\[
\xymatrix{ A \ar[r] \ar[d]_i & X \ar[d]^p \\
B \ar[r] & Y .}
\]
A lift exists in the diagram if the obstruction group $\dd ( \fib i
, \fib p )$ vanishes \cite[8.4]{obs}. By definition, $\fib i$
belongs to $\mm_{\geq n}$, and $\fib p$ belongs to $\mm_{ \leq n -1
}$. Hence the obstruction group vanishes because there are only
trivial maps in $\dd$ from objects in $\mm_{\geq n}$ to objects in
$\mm_{\leq n-1}$.

Now suppose that a map $i$ has the left lifting property with
respect to all co-$n$-fibrations. Lemma \ref{factor} allows us to
apply the retract argument and conclude that $i$ is a retract of an
$n$-cofibration. But $n$-cofibrations are preserved by retract by
Lemma \ref{cfretract}, so $i$ is an $n$-cofibration.

A similar argument shows that if $p$ has the right lifting property
with respect to all $n$-cofibrations, then $p$ is a
co-$n$-fibration.
\end{proof}

\begin{cor}
\label{cor:lift} The class of $n$-cofibrations is closed under
arbitrary cobase change. The class of co-$n$-fibrations is closed
under arbitrary base change.
\end{cor}

\begin{proof}
This follows immediately from Lemma \ref{lift} together with the
facts that cobase changes preserve left lifting properties and base
changes preserve right lifting properties.
\end{proof}

\begin{lem} \label{lem:n-cofibration-fibration}
Every acyclic fibration is a co-$n$-fibration. Every acyclic
cofibration is an $n$-cofibration. If $n \geq m$, then every
$n$-cofibration is an $m$-cofibration, and every co-$m$-fibration is
a co-$n$-fibration.
\end{lem}

\begin{proof} This follows from
part (2) and (3) of Lemma \ref{lem:W^n-prop}.
\end{proof}

\begin{lem}
\label{lem:susp-n-cofib} Let $f$ be a cofibration. Then $f$ is an
$n$-cofibration if and only if $f \wedge S^1$ is an
$(n+1)$-cofibration.
\end{lem}

\begin{proof}
First, $- \wedge S^1$ preserves cofibrations because the model
structure on $\mm$ is simplicial. Let $C$ be the cofiber of $f$;
this is also the homotopy cofiber of $f$ because $f$ is a
cofibration. Note that $C$ is cofibrant. Now $C \wedge S^1$ is the
cofiber (and also the homotopy cofiber) of $f \wedge S^1$ because
the model structure on $\mm$ is simplicial. Since $C$ is cofibrant,
$C \wedge S^1$ is homotopically correct and is a model for $\Sigma
C$ in $\dd$. Now $C$ belongs to $\mm_{\geq n+1}$ if and only if $C
\wedge S^1$ belongs to $\mm_{\geq n+2}$.
\end{proof}

We next show that $n$-cofibrations interact appropriately with the
simplicial structure. This will be needed to show that our later
constructions behave well simplicially.

\begin{pro} \label{prop:n-cofibration-pushout-product}
Suppose that $f\col A \map B$ is an $n$-cofibration and $i \col K
\map L$ is a cofibration of simplicial sets. Then the map
\[
g\col A \otimes L \amalg_{A \otimes K} B \otimes K
  \map B \otimes L
\]
is also an $n$-cofibration.
\end{pro}

\begin{proof}
The map $i$ is a transfinite composition of cobase changes of maps
of the form $\partial\Delta[j] \map \Delta[j]$. Therefore, the map
$g$ is a transfinite composition of cobase changes of maps of the
form
\[
A \otimes \Delta[j] \amalg_{A \otimes \partial\Delta[j]} B \otimes
\partial\Delta[j] \map B \otimes \Delta[j].
\]
Since $n$-cofibrations are characterized by a left lifting property
(see Lemma \ref{lift}), $n$-cofibrations are preserved by cobase
changes and transfinite compositions. Therefore, we may assume that
$i$ is the map $\partial\Delta[j] \map \Delta[j]$.

Since $\mm$ is a simplicial model category and $f$ is a cofibration,
$g$ is also a cofibration. We need only show that $g$ is an
$n$-equivalence.

Let $C$ be the cofiber of the $n$-cofibration $f$,
so $C$ belongs to $\mm_{\geq n+1}$.
Then the cofiber of $g$ is $C \wedge S^{j}$, where the simplicial
set $S^{j}$ is the sphere $\Delta[j] /
\partial\Delta[j]$ based at the image of $\partial\Delta[j]$.
We need to show that $C \wedge S^{j}$ also belongs to $\mm_{\geq n+1}$.
But $C \wedge S^j$ is a model for $\Sigma^j C$ in $\dd$
because $C$ is cofibrant, so $C \wedge S^j$ belongs to $\mm_{\geq n+1}$
because $\mm_{\geq n+1}$ is closed under $\Sigma$.
\end{proof}

Note that the reduced version of Proposition
\ref{prop:n-cofibration-pushout-product} also holds. Namely, if
$f\col A \map B$ is an $n$-cofibration and $i\col K \map L$ is a
cofibration of pointed simplicial sets, then the map
\[
A \wedge L \amalg_{A \wedge K} B \wedge K \map B \wedge L
\]
is also an $n$-cofibration. The proof is identical.

\begin{cor} \label{cor:n-Map}
Let $A \map B$ be an $n$-cofibration, and let $X \map Y$ be a
co-$n$-fibration. The map
\[
f\col \Map(B, X) \map \Map(A,X) \times_{\Map(A,Y)} \Map(B,Y)
\]
is an acyclic fibration of simplicial sets.
\end{cor}

\begin{proof}
This follows from the lifting property characterization of acyclic
fibrations, adjointness, and Proposition
\ref{prop:n-cofibration-pushout-product}.
\end{proof}

\subsection{Producing t-model categories}

We give some elementary results for constructing t-model structures.

\begin{lem} \label{t-structures}
Assume that $\dd$ is the homotopy category of a proper simplicial
stable model category $\mm$. Let $\dd_{\geq 0 }$ be a strictly full
subcategory of $\dd$ that is closed under $\Sigma$. Define
$\dd_{\geq n} $ to be $\Sigma^{n} \dd_{\geq 0}$. Let $W_n$ be
defined as in Definition \ref{defn:n-eq}, and set $F_n = \text{inj }
C \cap W_n$. Let $\dd_{\leq n-1 }$ be the full subcategory of $\dd$
whose objects are isomorphic to $\fib (g)$ for all $g$ in $F_n$. If
there is a functorial factorization of any map in $ \mm$ as a map in
$C\cap W_n$ followed by a map in $F_n$, then $\dd_{\geq 0}$,
$\dd_{\leq 0 }$ is a t-structure on $\dd$, and hence we get a
t-model structure on $\mm$.
\end{lem}

 \begin{proof}
  We verify that and $\dd_{\geq 0}$
$ \dd_{\leq 0}$ satisfy the three axioms of a t-structure on $\dd$
given in Definition \ref{defn:t}.
 Axiom 1 holds since $\dd_{\geq 0}$ is closed under $\Sigma$. The
factorization applied to
 $* \rarr X$ (or $X \rarr *$) gives a natural triangle fulfilling
axiom 2 for a t-structure.

 Now assume that $X \in \dd_{\geq
0}$ and that $Y \in \dd_{\leq - 1}$. We can assume that $X$ is
cofibrant. Factor $X \rarr
* $ into a cofibration $ g\col X \rarr Z$ followed by an acyclic
fibration $Z \rarr *$. We have that $g $ is in $C \cap W_0$ since
$g$ is a cofibration with homotopy cofiber in $\dd_{\geq 1} $. By
our assumption $Y$ is weakly equivalent to the pullback $Y'$ of a
fibration $p \col E \rarr B$ with fibrant target having the right
lifting property with respect to $g$. For any map $f \col X \rarr
Y'$ there are commutative squares
\[ \xymatrix{ X \ar[r] \ar[d]_g & Y' \ar[d] \ar[r] & E \ar[d]^p \\
Z \ar[r]^{\sim} & \ast \ar[r] & B.} \] We get that the left square
lifts by our assumptions. Hence any map $f \col X \rarr Y'$ factors
through a contractible object. Since $X$ is cofibrant and $Y'$ is
fibrant we get that $\dd ( X , Y') =0 $. Hence we conclude that $\dd
( X , Y) =0 $ whenever $ X \in \dd_{\geq 0}$ and $Y \in \dd_{\leq
-1}$.
\end{proof}

\begin{pro} \label{cofibrantlygenerated}
Let $\mm$ be a proper simplicial stable cofibrantly generated model
category with homotopy category $\dd$.
  Let $I$ be a set of
 generating cofibrations and let $J$ be a set of generating
acyclic cofibrations. Let $K_n $ be subsets of $J$ for $n \in \zz$.
 Let $C(n)$
be the class of retracts of relative $ I \cup K_n$-cell complexes.
Let $W_n$ be the corresponding class of $n$-equivalences defined as
the class of maps that is the composite of a map in $C(n)$ followed
by an acyclic fibration.

If $W_n$ is equal to $ \Sigma^{n} W_0$ for all $n$, then the
structure defined above is a t-model structure and $C(n) = C \cap
W_n$.
\end{pro}

\begin{proof} We have functorial factorization of any map as a
map in $C(n)$ followed by a map in $\inj C(n)$ \cite[10.5,
11.1.2]{hir}. The result follows from Lemma \ref{t-structures} by
letting $\dd_{\geq 0}$ be the full subcategory of $\dd$ consisting
of objects isomorphic to the homotopy fibers of maps in $W_0$.
 \end{proof}

\section{Review of pro-categories}
\label{sctn:pro-prelim}

We give a brief review of pro-categories. This section contains
mostly standard material on pro-categories \cite{SGA} \cite{AM}
\cite{EH}.

\begin{defn}
\label{defn:pro} For any category $\cc$, the category \mdfn{$\p
\cc$} has objects all cofiltering diagrams in $\cc$, and
\[
{\p \cc}(X,Y) = \lim_t \colim_s \,
   \cc (X_s, Y_t).
\]
Composition is defined in the natural way.
\end{defn}

A {\bf constant} pro-object is one indexed by the category with one
object and one (identity) map. Let $\mathbf{c}\col \cc \map \p \cc$
be the functor taking an object $X$ to the constant pro-object with
value $X$. Note that this functor makes $\cc$ a full subcategory of
$\p \cc$. The limit functor $\mdfn{\lim}\col \p \cc \map \cc$ is the
right adjoint of $c$.

A {\bf level map} $X \map Y$ is a pro-map that is given by a natural
transformation (so $X$ and $Y$ must have the same indexing
category); this is a very special kind of pro-map. Up to
pro-isomorphism, every map is a level map \cite[App.~3.2]{AM}.

Let $M$ be a collection of maps in $\cc$. A level map $g$ in $\p
\cc$ is a \mdfn{levelwise $M$-map} if each $g_s$ belongs to $M$. A
pro-map is an \mdfn{essentially levelwise $M$-map} if it is
isomorphic to a levelwise $M$-map.

We say that a level map is directed (resp., cofinite directed) if
its indexing category is a directed set (resp., cofinite directed
set). Recall that a directed indexing set $S$ is cofinite if for all
$s \in S$, the set $\{ t \in S \, | \, t<s \} $ is finite.

\begin{defn}
\label{defn:special} A map in $\p \cc$ is a \mdfn{special $M$-map}
if it is isomorphic to a cofinite directed levelwise map $f = \{ f_s
\}_{s \in S}$ with the property that for each $s \in S$, the map
\[
M_s f\col X_s \rarr \lim_{t <s } X_t \times_{ \lim_{t <s } Y_t} Y_s
\]
belongs to $M$.
\end{defn}

\subsection{Strict model structures}
\label{subsctn:strict}

If $\mm$ is a proper model category, then $\p \mm$ has a
\mdfn{strict model structure} \cite{EH} \cite{strict}. The
\mdfn{strict cofibrations} are the essentially levelwise
cofibrations, the \mdfn{strict weak equivalences} are the
essentially levelwise weak equivalences, and the \mdfn{strict
fibrations} are retracts of special fibrations (see Definition
\ref{defn:special}).

The functors $c$ and $\lim$ are a Quillen adjoint pair between $\mm$
and $\p \mm$. The right derived functor of $\lim$ is $\holim$
\cite[Rem.~4.2.11]{EH}. To see why this is true, recall that $R\lim
X$ is defined to be $\lim \hat{X}$, where $\hat{X}$ is a strict
fibrant replacement for $X$. Now $\hat{X} \map *$ is a special
fibration if and only if $\hat{X}$ is a ``Reedy fibrant'' diagram
\cite[Ch.~15]{hir}. This shows that $\lim \hat{X}$ is one of the
usual models for $\holim X$ \cite[XI]{bk}.

\begin{pro} \label{pro:mappingspace}
\label{prop:strict-map} Let $X$ be a cofibrant object of $\p \mm$.
Let $Y$ be any levelwise fibrant object of $\p \mm$ with strict
fibrant replacement $\hat{Y}$. Then the homotopically correct
mapping space $\Map(X,\hat{Y})$ in the strict model structure is
weakly equivalent to $\holim_t \colim_s \Map(X_s, Y_t)$.
\end{pro}

\begin{proof}
We may reindex $Y$ so that it is cofinite directed and still
levelwise fibrant \cite[Thm.~2.1.6]{EH}. Since $\Map(X, \hat{Y})$ is
homotopically correct, it doesn't matter which strict fibrant
replacement $\hat{Y}$ that we consider. Therefore, we may choose one
with particularly good properties. Use the method of
\cite[Lem.~4.7]{strict} to factor the map $Y \map *$ into a strict
acyclic cofibration $Y \map \hat{Y}$ followed by a special fibration
$\hat{Y} \map
*$. This particular construction gives that $Y \map \hat{Y}$
is a levelwise weak equivalence and that $\hat{Y}$ is levelwise
fibrant.

Define a new pro-space $Z$ by setting $Z_t = \colim_s \Map(X_s,
\hat{Y}_t)$. The map $\hat{Y}_t \map \lim_{u<t} \hat{Y}_u$ is a
fibration because $\hat{Y}$ is strict fibrant. Since finite limits
and directed colimits of simplicial sets commute, we get that the
map $Z_t \map \lim_{u<t} Z_u$ is a fibration, and $Z$ is a strict
fibrant pro-space. Therefore, the simplicial set $\Map(X, \hat{Y}) =
\lim_t Z_t$ is weakly equivalent to the simplicial set $\holim_t
Z_t$ because homotopy limit is the derived functor of limit.

The map $\colim_s \Map(X_s, Y_t) \map \colim_s \Map(X_s, \hat{Y}_t)$
is a weak equivalence because $Y_t \map \hat{Y}_t$ is a weak
equivalence between fibrant objects. Homotopy limits preserve
levelwise weak equivalences, so the map
\[
\holim_t \colim_s \Map(X_s, Y_t) \map \holim_t \colim_s \Map(X_s,
\hat{Y}_t)
\]
is a weak equivalence.
\end{proof}

\section{Model structures on pro-categories}
\label{sctn:pro-model}

The goal of this section is to construct a certain model structure
on $\p \mm$ when $\mm$ is a t-model structure. First, we make a
connection between t-model structures and \mdfn{filtered model
structures}. A filtered model structure is a highly technical
generalization of a model structure that is useful for producing
interesting model structures on pro-categories \cite{ffi}.

We denote by \mdfn{$F_n$} the class of co-$n$-fibrations in the
t-model structure $\mm$, and we write \mdfn{$C_n = C$} for the class
of cofibrations in $\mm$. Note that $C_n$ does not really depend on
$n$ and that it is {\em not} the class of $n$-cofibrations.
Recall from Definition \ref{defn:n-eq} that $W_n$ is the class of
$n$-equivalences.

\begin{pro} \label{filtered} Let $\mm$ be a t-model structure.
  Then $(W_n, C_n, F_n)$ is a
proper simplicial filtered model structure on $\mm$, where the
indexing set is $\zz$ with its usual ordering.
\end{pro}

\begin{proof}
We showed in part (2) of Lemma \ref{lem:W^n-prop} that $W_n$ is
contained in $W_m$ whenever $n \leq m$. The second half of part (2)
of Lemma \ref{lem:W^n-prop} implies that $F_m$ is contained in $F_n$
whenever $n \geq m$.

The class $\inj C$ of maps that have the right lifting property with
respect to $C$ is equal to the class of acyclic fibrations, while
the class $\proj F_n$ of maps that have the left lifting property
with respect to $F_n$ is equal to the class of $n$-cofibrations by
Lemma \ref{lift}. These observations are central to verification of
the axioms for a filtered model structure.

The axiom numbers below refer to \cite[Sec.~4]{ffi}. Axiom 4.2
follows from Lemma \ref{W^n-2outof3}. Axiom 4.3 follows from part
(5) of Lemma \ref{lem:W^n-prop}, Lemma \ref{cfretract}, and
Corollary \ref{cor:lift}. The first half of Axiom 4.4 follows from
our identification of $\inj C$ and by part (3) of Lemma
\ref{lem:W^n-prop}, while the second half is immediate from the
description of $\proj F_n$ as $C \cap W_n$. The first half of Axiom
4.5 is provided by factorizations into cofibrations followed by
acyclic fibrations, while the second half is Lemma \ref{factor}. For
Axiom 4.6, we can factor an $n$-equivalence into a cofibration
followed by an acyclic fibration; then the cofibration is
necessarily an $n$-cofibration. Axioms 4.9 and 4.10 are established
in Lemma \ref{proper}. The non-trivial part of Axiom 4.12 is
Proposition \ref{prop:n-cofibration-pushout-product}.
\end{proof}

\begin{defn}
\label{defn:pro-ms} A map in $\p \mm$ is an \mdfn{$\hh_*$-weak
equivalence} if it is an essentially levelwise $W_n$-equivalence for
all $n$. Let \mdfn{$F_{\infty}$} be the union $\cup_n F_n$. A map
in $\p \mm$ is an \mdfn{$\hh_*$-fibration} if it is a retract of a
special $F_{\infty}$-map.
\end{defn}

A justification for the terminology is given by Theorem
\ref{thm:pro-hh}, where we show that the $\hh_*$-weak equivalences
can be detected by the homology functors $\hh_n$. The notation
$F_{\infty}$ reflects the fact that $F_{n+1}$ contains $F_n$ for
all integers $n$.

The cofibrations in $\p \mm$ are the essentially levelwise
cofibrations. They are the same as the strict cofibrations (see
Section \ref{subsctn:strict}), so we do not need a new name for
them.

\begin{thm}
\label{thm:pro-ms} Let $\mm$ be a t-model structure. The essentially
levelwise cofibrations, $\hh_*$-weak equivalences, and
$\hh_*$-fibrations are a proper simplicial model structure on $\p
\mm$.
\end{thm}

This model structure on $\p \mm$ is called the \mdfn{$\hh_*$-model
structure}.

\begin{proof}
This follows immediately from Proposition \ref{filtered} and
\cite[Thms. 5.15,5.16]{ffi}.
\end{proof}

Theorem \ref{thm:pro-ms} applied to the Postnikov t-model structure
on a category of spectra gives the model structure on the category
of pro-spectra described in the introduction.
%%%We do not spell out
%%%the properties of the $\hh_*$-model structure on pro-chain
%%%complexes of $R$-modules for a ring $R$ with respect to the standard
%%%model structure.

\begin{rem} A t-structure is \dfn{constant} if $\dd_{\geq 0 } = \dd_{\geq 1}$.
Constant t-structures correspond to
triangulated localization functors. A localization functor is a
functor $L$ together with a natural transformation $\eta \col 1
\rarr L$ so that $ L \eta (X) = \eta (LX)$ and these maps are
isomorphisms for all $X \in \dd$.
%Let $L$ be an idempotent
%functor on $\dd$ that is left adjoint to the full inclusion
%functor $ L \dd \rarr \dd$. If $L$ is a triangulated functor,
%then there is a constant t-structure on $\dd$ so that $
%\dd^{\leq 0}$ is the full subcategory of colocal objects and $
%\dd^{\geq 1} $ is $ L \dd$.
The functor $\tau_{\leq 0} $ together with the natural
transformation $ \eta_0 \col 1 \rarr \tau_{\leq 0}$ is always a
localization functor. It is triangulated exactly when the
t-structure is constant.
A
constant t-model structure on a model category $\mm$ is a functorial
left Bousfield localization of $\mm$ with respect to the class of
maps $W_0$ \cite[3.3.1]{hir}.

The $\hh_*$-model structure associated to a constant t-model
structure on a category $\mm$ is the strict model structure on $\p
\mm$ obtained from the localized model structure on $\mm$.
\end{rem}

In order for the $\hh_*$-model structure to be useful, one needs a
better understanding of cofibrant objects and fibrant objects.
Moreover, an understanding of the $\hh_*$-acyclic cofibrations and
$\hh_*$-fibrations is also useful. We study these issues next.

Cofibrant objects are easy to describe. They are just essentially
levelwise cofibrant objects.

The following proposition gives useful criteria for detecting
$\hh_*$-acyclic cofibrations.

\begin{pro} \label{prop:hh-acyclic-cofibration}
Let $f$ be a map in $\p \mm$. The following are equivalent:
\begin{enumerate}
\item $f$ is an $\hh_*$-acyclic cofibration.
\item $f$ is an essentially levelwise $n$-cofibration for every $n$.
\item $f$ has the left lifting property with
respect to all constant pro-maps $cX \map cY$ in which $X \map Y$ is
a co-$m$-fibration for some $m$.
\end{enumerate}
\end{pro}

\begin{proof}
This follows from \cite[Prop.~4.11, 4.12]{ffi} and the lifting
property characterization of $n$-cofibrations given in Lemma
\ref{lift}.
\end{proof}

The $\hh_*$-fibrations are more difficult to describe. The
$\hh_*$-fibrations are strict fibrations since $F_{\infty}$ is
contained in the class of all fibrations. We take the strict
fibrations as our starting point and characterize the
$\hh_*$-fibrations among them.

\begin{lem} \label{fibrations} Let $p\col X \rarr Y$ be a special fibration
indexed on a cofinite directed set $S$. Then $p$ is a special
$F_{\infty}$-map if and only if for each $s \in S$ the map $p_s$ is
a co-$n$-fibration for some $n$.
\end{lem}

\begin{proof}
 We need to show that
each $M_s p \col X_s \rarr \lim_{t < s } X_t \times_{\lim_{t < s }
X_t} Y_s $ is in $F_{\infty}$ if and only if each $ p_s $ is in
$F_{\infty}$.

By induction, it suffices to prove that if both $M_t p$ and $p_t$
are in $F_{\infty}$ for all $t < s$, then $M_s p$ is in
$F_{\infty}$ if and only if $p_s$ is in $F_{\infty}$. We have a
pullback diagram
\[
\xymatrix{
\lim_{t < s } X_t \times_{\lim_{t < s } X_t} Y_s \ar[r] \ar[d] & Y_s \ar[d] \\
\lim_{t < s} X_t \ar[r] & \lim_{t <s }Y_t }
\]
in which the lower horizontal map is in $F_{\infty}$ since it is a
finite composition of base changes of the maps $M_t p$ for $t <s$
\cite[Lem.~2.3]{ffi}. Hence its base change $\lim_{t < s } X_t
\times_{\lim_{t < s } X_t} Y_s \rarr Y_s$ is also in $F_{\infty}$
by Corollary \ref{cor:lift}.

Now $p_s $ is the composition $ X_s \stackrel{M_s
p}{\longrightarrow} \lim_{t < s } X_t \times_{\lim_{t < s } X_t} Y_s
\longrightarrow Y_s$. If $M_s p$ belongs to $F_{\infty}$, then
$p_s$ is the composition of two maps in $F_{\infty}$ and is
therefore in $F_{\infty}$. On the other hand, if $p_s$ belongs to
$F_{\infty}$, then Lemma \ref{W^n-2outof3} implies that $M_s p$ is
a co-$n$-equivalence for some $n$. Since $M_s p$ was assumed to be a
fibration, this means that $M_s p$ is a co-$n$-fibration and hence
in $F_{\infty}$.
\end{proof}

We now give a characterization of $\hh_*$-fibrations. Let \mdfn{$\co
W_{\infty}$} be the union $\cup_n \co W_n$. The notation reminds us
that $\co W_{n+1}$ contains $\co W^n$ for all integers $n$.

\begin{pro} \label{characterizefibrations}
A map in $\p \mm$ is an $\hh_*$-fibration if and only if it is a
strict fibration and an essentially levelwise $\co W_{\infty}$-map.
\end{pro}

\begin{proof}
If $f \col X \rarr Y$ is a $\hh_*$-fibration, then $f$ is a strict
fibration since $F_{\infty}$ is contained in the class of all
fibrations. Since $f$ is a retract of a special $F_{\infty}$-map by
definition, it is a retract of an essentially levelwise
$F_{\infty}$-map because special $F_{\infty}$-maps are essentially
levelwise $F_{\infty}$-maps \cite[Lem.~5.14]{ffi}. Hence $f$ is
itself an essentially levelwise $F_{\infty}$-map since retracts
preserve essentially levelwise $F_{\infty}$-maps
\cite[Cor.~5.6]{lim}. Finally, just observe that $F_{\infty}$ is
contained in $\co W_{\infty}$.

For the other direction, we may assume that $f\col X \rarr Y$ is a
levelwise $\co W_{\infty}$-map indexed on a cofinite directed set.
Factor $f$ into a levelwise acyclic cofibration $i\col X \rarr Z$
followed by a special fibration $p\col Z \rarr Y$ using the method
of \cite[Lem.~4.7]{strict}. We have that $f$ is a retract of $p$.

By Lemma \ref{fibrations}, we just need to show that each map $p_s$
is a co-$n$-fibration for some $n$. Each $p_s$ is a fibration
because special fibrations are levelwise fibrations
\cite[Lem.~5.14]{ffi}. The map $p_s$ is a co-$n$-equivalence for
some $n$ by Lemma \ref{W^n-2outof3} (3) since $p_s i_s$ is $f_s$ and
$i_s$ is a weak equivalence.
\end{proof}

Finally, we are ready to identify the $\hh_*$-fibrant objects.

\begin{defn} \label{defn:bdbelow}
An object $X$ of $\mm$ is \mdfn{bounded above} if it belongs to
$\mm_{\leq n}$ for some $n$, and it is \mdfn{bounded below} if it
belongs to $\mm_{\geq n}$ for some $n$.
\end{defn}

\begin{pro}
\label{prop:fibrant} An object of $\p \mm$ is $\hh_*$-fibrant if and
only if it is strict fibrant and essentially levelwise bounded
above.
\end{pro}

\begin{proof}
This is immediate from Proposition \ref{characterizefibrations},
once we note that an object $X$ of $\mm$ is bounded above if and
only if the map $X \map *$ belongs to $\co W_{\infty}$.
\end{proof}

The following corollary simplifies the construction of
$\hh_*$-fibrant replacements.

\begin{cor}
\label{cor:fibrant} If $Y$ is an essentially levelwise bounded above
pro-object, then there is a strict fibrant replacement $\hat{Y}$ for
$Y$ such that $\hat{Y}$ is also a $\hh_*$-fibrant replacement for
$Y$.
\end{cor}

\begin{proof}
We may assume that $Y$ is levelwise bounded above and indexed on a
cofinite directed set. Factor the map $Y \map *$ into a levelwise
acyclic cofibration $i\col Y \map \hat{Y}$ followed by a special
fibration $p\col \hat{Y} \map *$ using the method of
\cite[Lem.~4.7]{strict}. Since $i$ is a levelwise weak equivalence,
it follows that $p$ is a levelwise $\co W_{\infty}$-map and thus a
levelwise $F_{\infty}$-map because it is a levelwise fibration
\cite[Lem.~5.14]{ffi}. Hence $p$ is a special $F_{\infty}$-map by
Proposition \ref{fibrations}, so $\hat{Y}$ is $\hh_*$-fibrant.
\end{proof}

\subsection{Quillen adjunctions}

In this section we give some conditions that guarantee that a
Quillen adjunction between two t-model structures $\mm$ and $\mm'$
gives a Quillen adjunction between the $\hh_*$-model structures on
$\p \mm$ and $\p \mm'$. We use this to show that the $\hh_*$-model
structure on $\p \mm$ is stable.

Recall that if $F\col \cc \map \cc'$ is a functor, then $F$ induces
another functor $\p \cc \map \p \cc'$ defined by applying $F$
levelwise. We will abuse notation and write $F$ also for this
functor. If $G$ is the right adjoint of $F$, then the induced
functor $G\col \p \cc' \map \p \cc$ is the right adjoint of $F\col
\p \cc \map \p \cc'$.

\begin{pro} \label{Quillenadjoint}
Let $\mm$ and $\mm'$ be two t-model categories and let $L \col \mm
\rarr \mm'$ be a left adjoint of $R \col \mm' \rarr \mm$. The
following are equivalent:
\begin{enumerate}
\item
The induced functors $L\col \p \mm \rarr \p \mm'$ and $R\col \p \mm'
\rarr \p \mm$ are a Quillen adjoint pair with respect to the
$\hh_*$-model structures on $\p \mm$ and $\p \mm'$.
\item
$L:\mm \map \mm'$ preserves cofibrations, and for every $n'$, there
is an $n$ such that $L$ takes $n$-cofibrations in $\mm$ to
$n'$-cofibrations in $\mm'$.
\item $R\col \mm' \map \mm$ preserves acyclic fibrations and also
preserves the class of maps that are co-$n$-fibrations for some $n$,
i.e., $R(F'_{\infty})$ is contained in $F_{\infty}$.
\end{enumerate}
\end{pro}

\begin{proof}
This follows from \cite[Lem.~3.7]{ffi}, \cite[Thm. 6.1]{ffi}, and
\cite[Prop.~6.2]{ffi}. Note that $\mm$ is a pointed model category,
and the classes of $n$-cofibrations are closed under retracts and
arbitrary small coproducts by Lemma \ref{lift}.
\end{proof}

\begin{pro} \label{Quillenequivalence}
Let $\mm$ and $\mm'$ be two t-model categories, and let $L \col \mm
\rarr \mm'$ be a left adjoint of $R \col \mm' \rarr \mm$ such that
the induced functors $L\col \p \mm \map \p \mm'$ and $R\col \p \mm'
\map \p \mm$ are a Quillen adjoint pair with respect to the
$\hh_*$-model structures on $\p \mm$ and $\p \mm'$. Assume also
that:
\begin{enumerate}
\item
For every $n'$, there is an $n$ such that if $X \rightarrow RY$ is
in $W_{n}$ with $X$ cofibrant in $\mm$ and $Y$ fibrant in $\mm'$,
then the adjoint map $LX \rightarrow Y$ is in $W'_{n'} $.
\item
For every $n$, there is an $n'$ such that if $L X \rightarrow Y$ is
in $W'_{n'}$ with $X$ cofibrant in $\mm$ and $Y$ fibrant in
$\mm'$, then the adjoint map $X \rightarrow R Y$ is in $W_n $.
\end{enumerate}
Then $L $ and $R$ induce a Quillen equivalence between the
$\hh_*$-model structures on $\p \mm$ and $\p \mm'$.
\end{pro}

\begin{proof}
This follows from Theorem \cite[Thm.~6.3]{ffi}.
\end{proof}

For any pro-object $X$, the suspension $X \wedge S^1$ is defined to
be the levelwise suspension of $X$, and the loops $\Map_*(S^1, X)$
is defined to be the levelwise loops of $X$.

The next theorem says that the $\hh_*$-model structure on $\p \mm$
is a {\em stable} model structure. In particular, the
$\hh_*$-homotopy category $\Ho( \p \mm )$ is a triangulated
category.

\begin{thm}
\label{thm:pi-Quillen-pair} The functors $- \wedge S^1$ and
$\Map_*(S^1,-)$ are a Quillen equivalence from the $\hh_*$-model
structure on $\p \mm$ to itself.
\end{thm}

\begin{proof}
On $\mm$, the functor $- \wedge S^1$ preserves cofibrations because
$\mm$ is a simplicial model category. By Lemma
\ref{lem:susp-n-cofib}, $- \wedge S^1$ takes $n$-cofibrations in
$\mm$ to $(n+1)$-cofibrations. Hence $- \wedge S^1$ is a left
Quillen adjoint by Proposition \ref{Quillenadjoint}. The goal of the
rest of the proof is to show that the conditions of Proposition
\ref{Quillenequivalence} are satisfied.

Let $g\col X \map \Map_*(S^1,Y)$ be any map in $\mm$ with $X$
cofibrant and $Y$ fibrant. Factor $g$ into a cofibration $i\col X
\map Z$ followed by an acyclic fibration $p\col Z \map
\Map_*(S^1,Y)$. The adjoint map $f\col X \wedge S^1 \map Y$ factors
as $X \wedge S^1 \map Z \wedge S^1 \map Y$, where the first map is
$i \wedge S^1$ and the second is adjoint to $p$. Note that the
second map is a weak equivalence because $\mm$ is a stable model
category.

Now $g$ is an $n$-equivalence if and only if $i$ is an
$n$-cofibration. Lemma \ref{lem:susp-n-cofib} implies that this
occurs if and only if $i \wedge S^1$ is an $(n+1)$-cofibration, and
this happens if and only if $f$ is an $(n+1)$-equivalence. Thus $g$
is an $n$-equivalence if and only if $f$ is an $(n+1)$-equivalence,
and both conditions of Proposition \ref{Quillenequivalence} are
satisfied.
\end{proof}

\section{Functorial towers of truncation functors}
\label{sctn:Postnikov}

Recall that in Lemma \ref{lem:rigid-truncation}, we showed that
there are rigid truncation functors $\tau_{\leq n}\col \mm \map
\mm_{\leq n}$ and $\tau_{\geq n}\col \mm \map \mm_{\geq n}$.
Although these functors are well-defined on the model category $\mm$
(not just on the homotopy category $\dd$), they have a major defect
as defined previously. Namely, there are no natural transformations
$\tau_{\leq n+1} \map \tau_{\leq n}$ and $\tau_{\geq n+1} \map
\tau_{\geq n}$. In this section, we will modify the definition of
the truncation functors so that these natural transformations do
exist.

Let $X$ be an object of $\mm$. We will define objects \mdfn{$T_{\leq
n} X$} for $n \geq 0$ inductively.

When $n = 0$, factor $X \map *$ into a $1$-cofibration $X \map
T_{\leq 0} X$ followed by a co-$1$-fibration $T_{\leq 0} X \map
*$. The object $T_{\leq 0} X$ belongs to $\mm_{\leq 0}$ by
definition of co-$1$-equivalences.

To define $T_{\leq 1} X$, factor the map $X \map T_{\leq 0} X$ into
a $2$-cofibration $X \map T_{\leq 1} X$ followed by a
co-$2$-fibration $T_{\leq 1} X \map T_{\leq 0} X$. Note that the
map $T_{\leq 1} X \map *$ is a composition of a co-$2$-fibration
with a co-$1$-fibration, so it is a co-$2$-equivalence by
Lemma \ref{W^n-2outof3}. Therefore, $T_{\leq 1} X$ belongs to
$\mm_{\leq 1}$ as desired.

Inductively, to construct $T_{\leq n} X$ for $n > 0$, assume that
$T_{\leq n-1} X$ and a natural map $X \map T_{\leq n-1} X$ have
already been constructed. Factor this map into an
$(n+1)$-cofibration $X \map T_{\leq n} X$ followed by a
co-$(n+1)$-fibration $T_{\leq n} X \map T_{\leq n-1} X$.

The construction is summarized by the tower
\[
\cdots \map T_{\leq 2} X \map T_{\leq 1} X \map T_{\leq 0} X.
\]

Note that we have not defined $T_{\leq n} X$ for $n < 0$. As far as
we know, it is not possible to define $T_{\leq n}$ for all $n$ so
that the desired natural transformations between these functors
exist.

Note also that $T_{\leq n}$ takes values in fibrant objects, and the
natural map $T_{\leq n} X \map T_{\leq n-1} X$ is a
co-$(n+1)$-fibration.

To define \mdfn{$T_{\geq n} X$} for $n\geq 1$, first recall that
maps in $\mm$ have functorial homotopy fibers. Then $T_{\geq n} X$
is defined to be the homotopy fiber of the map $X \map T_{\leq n-1}
X$. As before, we end up with a tower
\[
\cdots \map T_{\geq 3} X \map T_{\geq 2} X \map T_{\geq 1} X,
\]
and each $T_{\geq n} X$ belongs to $\mm_{\geq n}$.

The following lemma shows that $T_{\leq n} X$ and $T_{\geq n} X$
have the desired homotopy types.

\begin{lem}
\label{lem:rigid-truncate} The objects $T_{\leq n} X$ and
$\tau_{\leq n} X$ of $\mm$ are weakly equivalent. Similarly, the
objects $T_{\geq n} X$ and $\tau_{\geq n} X$ are weakly equivalent.
\end{lem}

\begin{proof}
There is a homotopy fiber sequence $T_{\geq n} X \map X \map
T_{\leq n-1} X$ such that $T_{\leq n-1} X$ belongs to $\mm_{\leq n-1}$
and $T_{\geq n} X$ belongs to $\mm_{\geq n}$. On $\dd$, $\tau_{\leq
n-1}$ and $\tau_{\geq n}$ are the unique functors with this property
(see Lemma \ref{truncation}). Therefore, $T_{\leq n-1}$ induces
$\tau_{\leq n-1}$ on $\dd$, which means that $T_{\leq n-1} X$ and
$\tau_{\leq n-1} X$ are weakly equivalent for all $X$. Similarly,
$T_{\geq n}$ induces $\tau_{\geq n}$ on $\dd$, so $T_{\geq n} X$ and
$\tau_{\geq n} X$ are weakly equivalent.
\end{proof}

\begin{lem}
\label{lem:Postnikov-tower} Let $Y$ be a pro-object indexed by a
cofiltered category $I$. Consider the pro-object $Z$ indexed on $I
\times \nn$ such that $Z_{s,n} = T_{\leq n} Y_s$. The natural map $Y
\map Z$ is an $\hh_*$-weak equivalence.
\end{lem}

\begin{proof}
For each $s$ and $n$, the map $Y_s \map Z_{s,n}$ has homotopy fiber
$T_{\geq n+1} Y_s$, so this map is an $(n+1)$-equivalence. This
shows that the map $Y \map Z$ is an essentially levelwise
$k$-equivalence for all $k$.
\end{proof}

The next result shows how the functors $T_{\leq n}$ are of
tremendous value in constructing $\hh_*$-fibrant replacements.

\begin{pro}
\label{pro:hh-fib-replacement} Let $Y$ be a pro-object indexed by a
cofiltered category $I$. Consider the pro-object $Z$ indexed on $I
\times \nn$ such that $Z_{s,n} = T_{\leq n} Y_s$. A strict fibrant
replacement for $Z$ is an $\hh_*$-fibrant replacement for $Y$.
\end{pro}

\begin{proof}
In order to construct an $\hh_*$-fibrant replacement for $Y$, Lemma
\ref{lem:Postnikov-tower} says that we may construct an
$\hh_*$-fibrant replacement for $Z$ instead. Finally, Corollary
\ref{cor:fibrant} says that a strict fibrant replacement for $Z$ is
the desired $\hh_*$-fibrant replacement.
\end{proof}

\section{Homotopy classes of maps of pro-spectra}
\label{sctn:homotopy-class}

We continue to work in a t-model structure $\mm$. Recall that $\dd$
is the homotopy category of $\mm$. Let \mdfn{$\ee$} be the
$\hh_*$-homotopy category of $\p \mm$.

The mapping space $\Map(X,Y)$ is related to homotopy classes in the
following way \cite[6.1.2]{hov}. For every cofibrant $X$ and
$\hh_*$-fibrant $Y$, $\ee (X, Y )$ is isomorphic to $\pi_{0} \Map(X,
Y)$.

%\begin{lem} \label{lem:mappingspace}
%Let $X$ be a pro-object in $\p \mm$, and let $Y$ be an
%essentially bounded below pro-object in $\p \mm$. Then the
%homotopically correct mapping spaces from $X$ to $Y$ in the
%strict and in the $\hh_*$-model structure are weakly
%equivalent.
%\end{lem}
%
%\begin{proof}
%This follows from Corollary \ref{cor:fibrant}.
%\end{proof}

\begin{lem} \label{cQuillen}
When $\p \mm$ is equipped with the $\hh_*$-model structure, the
constant pro-object functor $c \col \mm \rarr \p \mm$ and the limit
functor $\lim\col \p \mm \map \mm$ are a Quillen adjoint pair.
\end{lem}

\begin{proof}
Note that $c$ preserves cofibrations and acyclic cofibrations.
\end{proof}

\begin{pro}
The right derived functor $R \lim$ of $lim\col \p \mm \map \mm$ is
given by $R\lim Y = \holim_{s ,n} T_{\leq n} Y_{s}$.
\end{pro}

\begin{proof}
We may assume that $Y$ is indexed by a cofinite directed set $I$.
Let $Z$ be the pro-object indexed by $I \times \nn$ such that
$Z_{s,n}$ equals $T_{\leq n} Y_s$. Recall from Lemma
\ref{lem:Postnikov-tower} that the natural map $Y \map Z$ is an
$\hh_*$-weak equivalence.

Let $\hat{Z}$ be a strict fibrant replacement for $Z$. Corollary
\ref{cor:fibrant} says that $\hat{Z}$ is an $\hh_*$-fibrant
replacement for $Y$, so $R\lim Y$ is equal to $\lim \hat{Z}$. As
observed in Section \ref{subsctn:strict}, $\lim \hat{Z}$ is the same
as $\holim Z$.
\end{proof}

\begin{cor} \label{prop:constant}
There is a natural isomorphism
\[
\ee ( c X , Y ) \cong \dd ( X , \holim_{t, n} T_{\leq n} Y_t )
\]
for all $X$ in $\mm$ and all $Y$ in $\p \mm$. There is a natural
isomorphism
 \[
\ee ( cX , cY ) \cong \dd ( X , \holim_{ n \rarr \infty} T_{\leq
n} Y)
\]
for all $X$ and $Y$ in $\mm$.
\end{cor}

Consequently, if $ Y \rarr \holim_{n \rarr \infty} T_{\leq n} Y$
is a weak equivalence for all $Y$ in $\mm$, then the homotopy
category of $\mm$ embeds into the $\hh_*$-homotopy category on $\p
\mm$.

\begin{pro} \label{cor:homset}
Let $X$ and $Y$ be objects in $\p \mm$. Let $\hat{X}$ be a cofibrant
replacement of $X$.
% and let $\hat{Y}$ be a strict
%fibrant replacement of $\{ T^{\geq n} Y_t \}$ in $\p \mm$.
The homotopically correct mapping space of maps from $ X$ to $Y$ in
the $\hh_*$-model structure is weakly equivalent to the space
$\holim_{t ,n } \colim_s \Map(\hat{X}_s, T_{\leq n} Y_t)$.
\end{pro}
\begin{proof} This follows by
Propositions \ref{prop:strict-map} and \ref{pro:hh-fib-replacement}.
\end{proof}

Let $ X$ and $Y$ be objects in $ \p \mm$. In general the group $\ee
(X ,Y)$ is quite different from $\p \dd ( X , Y )$. There is not
even a canonical map from one to the other. The next Lemma says that
under some strong conditions the homsets in $\ee$ and $ \p \dd$
agree. We choose to use conceptual proofs rather than the higher
derived limit spectral sequence relating $\ee (X ,Y)$ to higher
derived limits of the inverse system $\{ colim_s \dd ( X_s , T_{\leq
n} Y_t) \}_{t,n} $ of abelian groups.

\begin{lem}
\label{heartmorphisms} Let $X$ and $Y$ be two pro-objects such that
$X_s$ is in $\mm_{\geq n}$ for all $s$ and $Y_t$ is in $\mm_{\leq n
}$ for all $t$. Then $\ee( X ,Y )$ is isomorphic to $\lim_t \colim_s
\dd ( X_s , Y_t ).$
\end{lem}

\begin{proof}
We may assume that $X$ is cofibrant. By taking a levelwise fibrant
replacement, we may assume that each $Y_s$ is fibrant. Corollary
\ref{cor:fibrant} and Proposition \ref{pro:mappingspace} imply that
the homotopically correct mapping space of maps from $X$ to $Y$ is
$\holim_t \colim_s \Map(X_s, Y_t)$, and we want to compute $\pi_0$
of this space.

For $k \geq 1$, the only map $\Sigma^k X_s \map Y_t$ in $\dd$ is the
trivial map because $\Sigma^k X_s$ belongs to $\dd_{\geq n+k}$ while
$Y_t$ belongs to $\dd_{\leq n}$. Therefore, $\pi_0 \Map(\Sigma^k
X_s, Y_t) = \pi_k \Map(X_s, Y_t)$ is trivial. We have just shown
that $\Map(X_s, Y_t)$ is a homotopy-discrete space.

Filtered colimits preserve homotopy-discrete spaces; moreover, they
commute with $\pi_0$. Similarly, homotopy limits preserve
homotopy-discrete spaces and respect $\pi_0$ in the sense that the
set of components of a homotopy limit is the ordinary limit of the
sets of components of each space. In our situation, this implies
that $\pi_0 \holim_t \colim_s \Map(X_s, Y_t)$ is equal to $\lim_t
\colim_s \pi_0 \Map(X_s, Y_t)$, which is the desired result.
\end{proof}

\begin{cor} \label{0mappingspace}
Let $X$ and $Y$ be two pro-objects such that $X_s$ is in $\mm_{\geq
n}$ for all $s$ and $Y_t$ is in $\mm_{\leq n-1}$ for all $t$. Then
$\ee( X ,Y )$ is trivial.
\end{cor}

\begin{proof}
This follows from Lemma \ref{heartmorphisms} and part (3) of
Definition \ref{defn:t}.
\end{proof}

\begin{cor} \label{cor:XcY}
If $Y$ is a bounded above object in $\mm$ and $X$ is any object in
$\p \mm$, then $\ee( X , cY)$ is isomorphic to $\colim_s \dd ( X_s ,
Y )$.
\end{cor}

\begin{proof}
The proof is nearly the same as the proof of Lemma
\ref{heartmorphisms}. We need to compute $\pi_0 \colim_s \Map(X_s,
Y)$. We just need to observe that $\pi_0$ commutes with filtered
colimits.
\end{proof}

Since $\p \mm$ is a model category, one may consider homotopy limits
internal to $\p \mm$. In other words, given a diagram of
pro-objects, one can form the homotopy limit of this diagram and
obtain another pro-object. We will need the following basic result
about homotopy limits of countable towers later when we discuss
convergence of spectral sequences.

\begin{lem} \label{lim1} Let $\mm$ be a simplicial model
category, and let $\dd$ be its homotopy category. Let $X$ belong to
$ \mm$, and let
\[
\cdots \map Y^2 \map Y^1 \map Y^0
\]
be a countable tower in $ \mm$. There is a natural short exact
sequence
\[
\lim_k^1 \dd ( \Sigma X , Y^k ) \rarr \dd ( X , \holim_k Y^k ) \rarr
\lim_k \dd (X , Y^k ). \]
\end{lem}

\begin{proof}
We may assume that $X$ is cofibrant, that each $Y^k$ is
$\hh_*$-fibrant, and that each map $Y^k \map Y^{k-1}$ is an
$\hh_*$-fibration. We have that $ \Map ( X , \lim_k Y_k ) $ is
isomorphic to $ \lim_k \Map ( X , Y_k ) $ since $\Map (X , - )$ is
right adjoint to tensoring with $X$. Since $ - \otimes X$ sends
acyclic cofibrations of simplicial sets to acyclic cofibrations in
$\mm$, we get that the tower $\Map ( X , Y_k ) $ is a tower of
fibrations between fibrant simplicial sets.

Hence $ \Map ( X , \lim_k Y_k ) $ is equivalent to $\holim_k \Map(X,
Y^k)$. The claim now follows by the $\lim^1$ short exact sequence
for simplicial sets \cite[IX.3.1]{bk}.
\end{proof}

\section{t-model structure for pro-categories}
\label{sctn:t-structureonprocategoreis} We now define a t-structure
on the $\hh_*$-homotopy category $\ee$ of $\p \mm$.

\begin{defn}
\label{defn:pro-t} Let \mdfn{$(\p \mm)_{\leq 0}$} be the full
subcategory of $\p \mm$ on all objects that are $\hh_*$-weakly
equivalent to a pro-object $X$ such that each $X_s$ belongs to
$\mm_{\leq 0}$. Let \mdfn{$(\p \mm)_{\geq 0}$} be the full
subcategory of $\p \mm$ on all objects that are $\hh_*$-weakly
equivalent to a pro-object $X$ such that each $X_s$ belongs to
$\mm_{\geq 0}$.
\end{defn}

We define \mdfn{$(\p \mm)_{\leq n}$} and \mdfn{$(\p \mm)_{\geq n}$}
to be the subcategories $\Sigma^{n} (\p \mm)_{\leq 0}$ and
$\Sigma^{n} (\p \mm)_{\geq 0}$ respectively. Recall that $\Sigma$
here refers to the levelwise suspension functor on pro-objects.

\begin{lem}
\label{lem:pro-leq-n} The subcategory $(\p \mm)_{\leq n}$ is the
full subcategory of $\p \mm$ on all objects that are $\hh_*$-weakly
equivalent to a pro-object $X$ such that each $X_s$ belongs to
$\mm_{\leq n}$. Similarly, the subcategory $(\p \mm)_{\geq n}$ is
the full subcategory of $\p \mm$ on all objects that are
$\hh_*$-weakly equivalent to a pro-object $X$ such that each $X_s$
belongs to $\mm_{\geq n}$.
\end{lem}

\begin{proof}
We prove the first claim. The proof of the second claim is similar.

First suppose that $X$ is a pro-object such that each $X_s$ belongs
to $\mm_{\leq n}$. Since $\Sigma^{-n}$ takes $\mm_{\leq n}$ to
$\mm_{\leq 0}$, it follows that $\Sigma^{-n} X$ belongs to $\mm_{\leq
0}$ levelwise. Therefore $\Sigma^{-n} X$ belongs to $(\p \mm)_{\leq
0}$, and $X$ belongs to $\Sigma^{n} (\p \mm)_{\leq 0}$.

Now suppose that $Y$ belongs to $\Sigma^{n} (\p \mm)_{\leq 0}$. It
follows that $\Sigma^{-n} Y$ belongs to $(\p \mm)_{\leq 0}$, so it is
$\hh_*$-weakly equivalent to a pro-object $X$ such that each $X_s$
belongs to $\mm_{\leq 0}$. Note that $\Sigma^{n} X$ belongs to
$\mm_{\leq n}$ levelwise. But $\Sigma^{n} X$ is $\hh_*$-weakly
equivalent to $Y$. This is the desired result.
\end{proof}

\begin{defn}
Let \mdfn{$\ee_{\leq n}$} be the full subcategory of $\ee$ on all
objects whose $\hh_*$-weak homotopy types belong to $(\p \mm)_{\leq
n}$. Let \mdfn{$\ee_{\geq n}$} be the full subcategory of $\ee$ on
all objects whose $\hh_*$-weak homotopy types belong to $(\p
\mm)_{\geq n}$.
\end{defn}

\begin{pro} \label{pro:tprostructure}
The classes $\ee_{\geq 0} $ and $\ee_{\leq 0}$ are a t-structure on
the $\hh_*$-homotopy category $\ee$ of $\p \mm$. Moreover, $\cap_n
\ee_{ \geq n}$ consists only of contractible objects.
\end{pro}

\begin{proof}
We verify the axioms in Definition \ref{defn:t}.

For part (1), suppose that $X$ belongs to $\ee_{\geq 0}$. We may
assume that each $X_s$ belongs to $\mm_{\geq 0}$. Now each $\Sigma
X_s$ belongs to $\mm_{\geq 0}$, so $\Sigma X$ belongs to $\mm_{\geq
0}$ levelwise. Thus $\Sigma X$ lies in $\ee_{\geq 0}$. To show that
$\ee_{\leq 0}$ is closed under $\Sigma^{-1}$, use the dual argument.

Lemma \ref{lem:pro-leq-n} implies that $\Sigma^{-1} \ee_{\leq 0}$ is
the full subcategory of $\ee$ on all objects that are $\hh_*$-weakly
equivalent to an object $X$ such that each $X_s$ belongs to
$\mm_{\leq -1}$. This will be needed in parts (2) and (3) below.

For part (2), let $X$ be any object in $\p \mm$. Apply the
truncation functors $\tau_{\geq 0}$ and $\tau_{\leq -1}$ to obtain a
levelwise homotopy cofiber sequence $\tau_{\geq 0} X \map X \map
\tau_{\leq -1} X$. Finally, observe that levelwise homotopy cofiber
sequences are homotopy cofiber sequences in the $\hh_*$-model
structure because levelwise cofibrations are cofibrations. Note that
we need Lemma \ref{lem:pro-leq-n} to conclude that $\tau_{\leq -1} X$
belongs to $\ee_{\leq -1}$.

Part (3) is Corollary \ref{0mappingspace}, again using Lemma
\ref{lem:pro-leq-n} to identify $\ee_{\leq -1}$.

For the last claim, suppose that $X$ belongs to $\cap_n \ee_{\geq
n}$. Fix a value of $n$. Then we may assume that each $X_s$ belongs
to $\mm_{\geq n}$, so the map $X \map *$ is a levelwise
$n$-equivalence. Thus $X \map *$ is an $\hh_*$-weak equivalence, so
$X$ is contractible.
\end{proof}

The subcategory $\cap_n \ee_{\geq n}$ contains only contractible
objects even if $\cap_n \dd_{\geq n}$ contains non-contractible
objects. On the other hand, $\cap_n \ee_{\leq n}$ contains only
contractible objects if and only if $\dd = \dd_{\geq 0}$.

\begin{lem} \label{rem:constant} If all the objects of $\cap_n
\ee_{\leq n}$ are contractible, then $\dd$ is equal to $ \dd_{\geq
0}$.
\end{lem}

\begin{proof}
Assume that there are noncontractible elements $ X_m \in \mm_{\leq
m}$ in each degree $m$. Define a pro-object $\{ Y_n \}$ by letting
$Y_n = \coprod_{ m \leq n} X_m $ and letting the map $ Y_{n-1} \rarr
Y_{n}$ be the canonical map ($\mm$ has a zero object). We have that
$\{ Y_n \}$ is in $\cap_n ( \p \mm)_{\leq n}$, but $\{ Y_n \}$ is
noncontractible in $\p \mm$: If there is a weak equivalence between
$\{ Y_n \}$ and $*$ in $\ee$, then for every $n$ there are integers
$n'$ and $m$ such that in the homotopy category $\dd$ of $\mm$ the
map $ (Y_{n'})_{\leq m} \rarr (Y_{n})_{\leq m}$ is the zero map. This
gives a contradiction since $(Y_{n'})_{\leq m}$ is not contractible
in $ \dd$ for any $m$.
\end{proof}

We can now identify the heart $\hh(\ee)$ of the t-structure from
Proposition \ref{pro:tprostructure}.

\begin{lem}
\label{lem:pro-heart} The category $\ee_{\leq 0} \cap \ee_{\geq 0}$
is the $\hh_*$-homotopy category of the subcategory $\p (\mm_{\leq
0} \cap \mm_{\geq 0})$ of $\p \mm$.
\end{lem}

\begin{proof}
Let $X$ be an object of $\ee_{\leq 0} \cap \ee_{\geq 0}$. We need to
show that $X$ is $\hh_*$-weakly equivalent to an object $Y$ of $\p
(\mm_{\leq 0} \cap \mm_{\geq 0})$.

We may assume that $X$ belongs to $\mm_{\geq 0}$ levelwise. If we
apply the functors $\tau_{\geq 1}$ and $\tau_{\leq 0}$ to $X$
levelwise, we obtain a levelwise homotopy cofiber sequence
$\tau_{\geq 1} X \map X \map \tau_{\leq 0} X$. This is a homotopy
cofiber sequence in $\p \mm$ because cofibrations are defined to be
levelwise cofibrations. Therefore,
\[
\tau_{\geq 1} X \map X \map \tau_{\leq 0} X \map \Sigma \tau_{\geq
1} X
\]
is a distinguished triangle in $\ee$. The object $\tau_{\geq 1} X$
belongs to both $\ee_{\geq 1}$ and to $\ee_{\leq 0}$, so it is
contractible. This means that the map $X \map \tau_{\leq 0} X$ is an
$\hh_*$-weak equivalence.

Finally, $\tau_{\leq 0} X$ belongs to $\mm_{\leq 0}$ levelwise by
definition of $\tau_{\leq 0}$, and it also belongs to $\mm_{\geq 0}$
levelwise because $X$ belongs to $\mm_{\geq 0}$ levelwise. Thus,
$\tau_{\leq 0} X$ is the desired pro-object $Y$.
\end{proof}

As in Definition \ref{defn:n-eq}, we define an
\mdfn{$n$-equivalence} (resp., \mdfn{co-$n$-equivalence}) in $\p
\mm$ to be a map whose homotopy fiber belongs to $(\p \mm)_{\geq n}$
(resp., homotopy cofiber belongs to $\p \mm)_{\leq n}$). By
definition, the $n$-equivalences in $\p \mm$ are different than the
levelwise $n$-equivalences. A similar warning applies to
co-$n$-equivalences. However, we will show below in Lemma
\ref{lem:pro-n-eq} that actually they coincide.

Unfortunately, we cannot conclude that $\p \mm$ has a t-model
structure. Although we can factor any map in $\p \mm$ into an
$n$-equivalence followed by a co-$n$-equivalence, it does not seem
to be possible to make this factorization functorial. Absence of
functorial factorizations is a general problem with pro-categories.
However, we will prove a slightly weaker result below in Proposition
\ref{pro:pro-t-model}.

\begin{lem}
\label{lem:level-n-eq} If $f$ is an essentially levelwise
$n$-equivalence in $\p \mm$, then $f$ is an $n$-equivalence in $\p
\mm$. If $f$ is an essentially levelwise co-$n$-equivalence in $\p
\mm$, then $f$ is a co-$n$-equivalence in $\p \mm$.
\end{lem}

\begin{proof}
Let $f$ be a levelwise $n$-equivalence. Homotopy cofibers of
pro-maps can be computed levelwise because cofibrations are defined
levelwise. It follows that $\hocofib f$ belongs to $\mm_{\geq n+1}$
levelwise. Lemma \ref{lem:pro-leq-n} says that $\hocofib f$ belongs
to $(\p \mm)_{\geq n+1}$. By definition, $f$ is an $n$-equivalence.

A similar argument proves the second claim.
\end{proof}

\begin{pro}
\label{pro:pro-t-model} The $\hh_*$-model structure on $\p \mm$ and
the t-structure on $\ee$ of Definition \ref{defn:pro-t} are a {\em
non-functorial} t-model structure on $\p \mm$ in the sense that all
the axioms of a t-model structure are satisfied except that the
factorizations in the model structure and the factorizations into
$n$-equivalences followed by co-$n$-equivalences might not
necessarily be functorial.
\end{pro}

\begin{proof}
We showed in Theorems \ref{thm:pro-ms} and \ref{thm:pi-Quillen-pair}
that the $\hh_*$-model structure is simplicial, proper, and stable.
We showed in Proposition \ref{pro:tprostructure} that Definition
\ref{defn:pro-t} is a t-structure on $\ee$.

It remains only to produce (non-functorial) factorizations into
$n$-equivalences followed by co-$n$-equivalences. Let $f \col X \map
Y$ be any map in $\p \mm$, which we may assume is a levelwise map.
Using that the t-model structure on $\mm$ has functorial
factorizations, we may factor $f$ into a levelwise $n$-equivalence
$g\col X \map Z$ followed by a levelwise co-$n$-equivalence $h\col Z
\map Y$. Finally, Lemma \ref{lem:level-n-eq} implies that $g$ is an
$n$-equivalence in $\p \mm$, and $h$ is a co-$n$-equivalence in $\p
\mm$.
\end{proof}

\begin{lem}
\label{lem:pro-n-eq} A map in $\p \mm$ is an $n$-equivalence if and
only if it is an essentially levelwise $n$-equivalence. A map in $\p
\mm$ is a co-$n$-equivalence if and only if it is an essentially
levelwise co-$n$-equivalence.
\end{lem}

\begin{proof}
We prove the first claim. The proof of the second claim is dual.

One direction was already proved in Lemma \ref{lem:level-n-eq}. For
the other direction, suppose that $f\col X \map Y$ is an
$n$-equivalence in $\p \mm$. We may assume that $f$ is a directed
cofinite level map. Factor $f$ into a levelwise cofibration $i\col X
\map Z$ followed by a special acyclic fibration $p\col Z \map Y$.
The map $p$ is an $\hh_*$-weak equivalence, so $i$ is also an
$n$-equivalence in $\p \mm$. Moreover, the map $p$ is a levelwise
weak equivalence. The class of essentially levelwise
$n$-equivalences is closed under composition (the proof of
\cite[Lem.~3.5]{strict} applies), so it suffices to show that $i$ is
an essentially levelwise $n$-equivalence.

We know that $i$ is an $n$-cofibration in $\p \mm$, so it has the
left lifting property with respect to all co-$n$-fibrations. Factor
the map $i$ into a levelwise $n$-cofibration $j\col X \map W$
followed by a special co-$n$-fibration $q\col W \map Z$. The map $q$
is a co-$n$-fibration in $\p \mm$, so $i$ has the left lifting
property with respect to $q$ by Lemma \ref{lift}. Thus, the retract
argument shows that $i$ is a retract of $j$. But essentially
levelwise $n$-cofibrations are closed under retract
\cite[Cor.~5.6]{lim}, so $i$ is an essentially levelwise
$n$-cofibration and thus an essentially levelwise $n$-equivalence.
\end{proof}

\begin{lem}
\label{lem:hh-strict} The strict homotopy category of $\p (\mm_{\leq
0} \cap \mm_{\geq 0})$ is the same as the $\hh_*$-homotopy category
of $\p (\mm_{\leq 0} \cap \mm_{\geq 0})$.
\end{lem}

\begin{proof}
In order to compute strict weak homotopy classes from $X$ to $Y$,
one needs to take a strict cofibrant replacement $\tilde{X}$ of $X$
and a strict fibrant replacement $\hat{Y}$ of $Y$. But $\tilde{X}$
is also an $\hh_*$-cofibrant replacement for $X$, and Corollary
\ref{cor:fibrant} says that $\hat{Y}$ is an $\hh_*$-fibrant
replacement for $Y$.
\end{proof}

We would like a description of the heart of the t-structure on
$\ee$. We have not been able to identify the heart in complete
generality. However, in the primary applications to chain complexes
or to spectra, we can identify it using Proposition
\ref{pro:rigid-heart}

\begin{pro}
\label{pro:rigid-heart} Suppose that there is a ``rigidification''
functor $K\col \hh(\dd) \map \mm_{\leq 0} \cap \mm_{\geq 0}$ such
that the composition $\hh(\dd) \map \mm_{\leq 0} \cap \mm_{\geq 0}
\map \hh(\dd)$ is the identity. Then the heart $\hh(\ee)$ of the
$\hh_*$-homotopy category on $\p \mm$ is equivalent to the category
$\p \hh(\dd)$.
\end{pro}

For spectra, the functor $K$ takes an abelian group $A$ to a
functorial model for the Eilenberg-Mac\,Lane spectrum $HA$. For
chain complexes, $K$ takes an $R$-module $A$ to the chain complex
with value $A$ concentrated in degree $0$.

\begin{proof}
The functor $K$ extends to a levelwise functor $\p \hh(\dd) \map \p
(\mm_{\leq 0} \cap \mm_{\geq 0})$. This gives us a functor $F\col \p
\hh(\dd) \map \hh(\ee)$ after composition with the usual quotient
functor (because the quotient functor takes $\p (\mm_{\leq 0} \cap
\mm_{\geq 0})$ into both $\ee_{\leq 0}$ and $\ee_{\geq 0}$).

On the other hand, the quotient functor $\mm_{\leq 0} \cap \mm_{\geq
0} \map \hh(\dd)$ extends to a levelwise functor $\p (\mm_{\leq 0}
\cap \mm_{\geq 0}) \map \p \hh(\dd)$. This functor takes levelwise
weak equivalences to (levelwise) isomorphisms, so the functor
factors through the strict homotopy category of $\p (\mm_{\leq 0}
\cap \mm_{\geq 0})$. Lemma \ref{lem:hh-strict} implies that the
functor also factors through the $\hh_*$-homotopy category of $\p
(\mm_{\leq 0} \cap \mm_{\geq 0})$, which is the same as $\hh(\ee)$
by Lemma \ref{lem:pro-heart}. Thus we obtain a functor $G\col
\hh(\ee) \map \p \hh(\dd)$.

It remains to show that $F$ and $G$ are inverse equivalences. The
composition $FG$ is the identity because of the original assumption
on $K$. On the other hand, for every pro-object $X$, $GFX$ is
levelwise weakly equivalent to $X$. Thus $GF$ is isomorphic to the
identity.
\end{proof}

Without a rigidification functor, the most we can say is stated in
the following lemma.

\begin{lem}
\label{lem:full-faithful} The functor $G\col \hh(\ee) \map \p
\hh(\dd)$ from the proof of Proposition \ref{pro:rigid-heart} is
fully faithful.
\end{lem}

\begin{proof}
See Lemma \ref{heartmorphisms}.
\end{proof}

As promised in Section \ref{sctn:pro-model}, we now recharacterize
$\hh_*$-weak equivalences in terms of the pro-homology functors
$\hh_n$.

\begin{thm}
\label{thm:pro-hh} A map $f\col X \map Y$ in $\p \mm$ is an
$\hh_*$-weak equivalence if and only if it is an essentially
levelwise $m$-equivalence for some $m$ and $\hh_n(f)$ is an
isomorphism in $\p \hh(\dd)$ for all $n$.
\end{thm}

\begin{proof}
This is an application of Theorem \ref{thm:Whitehead} to the
$\hh_*$-model structure. The hypothesis of that theorem is proved at
the end of Proposition \ref{pro:tprostructure}. We also need Lemma
\ref{lem:pro-n-eq} to identify the $m$-equivalences in $\p \mm$.
Finally, we need Lemma \ref{lem:full-faithful} to recognize that a
map $g$ is an isomorphism in $\hh(\ee)$ if and only if $G(g)$ is an
isomorphism in $\p \hh(\dd)$ (where $G$ is the functor of Lemma
\ref{lem:full-faithful}).
\end{proof}

\section{The Atiyah-Hirzebruch spectral sequence}
\label{ah}

In this section we construct a spectral sequence for computing in
the homotopy category of a t-model structure. We will also
specialize this construction to the case of $\hh_*$-model structures
on pro-categories. Applied to the homotopy category of spectra with
the Postnikov t-structure we recover the Atiyah-Hirzebruch spectral
sequence for spectra.

 Recall from Definition \ref{homology}
 that $\tau_{\leq q} \tau_{\geq q} Y$ is isomorphic to $
\Sigma^{q} \hh_{q} (Y) $, where $\hh_q (Y)$ is the $q$-th homology
of $Y$ with values in the heart of the t-structure. Also recall from
Definition \ref{Ecohomology} that $H^{-p} (X ; E ) $ is the $-p$-th
cohomology of $X$ with $E$-coefficients, where $E$ belongs to the
heart $\hh (\dd)$.

For brevity we write $[X,Y]_n$ instead of $\dd ( X , \Sigma^{-n} Y)$.

\begin{thm} \label{specseq}
For any $X$ and $Y$ in a t-model category $\mm$, there is a spectral
sequence with
\[ E_{p,q}^2 = H^{-p} ( X ; \hh_{q} ( Y) ).
\]
The spectral sequence conditionally converges to $[ X , Y ]_{p +q }$
if $\displaystyle\holim_{q\rarr \infty} T_{\geq q} Y$ is
contractible and if $X$ is bounded below.
\end{thm}

The construction of the spectral sequence is standard (for example,
see \cite[Sec. 12]{B} or \cite[App.~B]{gm}). Conditional convergence
of spectral sequences is defined in \cite[Defn.~5.10]{B}.

To set up the spectral sequence, we only use the t-structure on the
homotopy category, but we use homotopy theory to state the
convergence criterion. We are using the functors $T_{\geq q}$ rather
than the functors $\tau_{\geq q}$ so that we obtain an actual tower
\[
\cdots \map T_{\geq q+1} Y \map T_{\geq q} Y \map T_{\geq q-1} Y
\map \cdots
\]
whose homotopy limit we can consider. Recall that there are not
necessarily maps $\tau_{\geq q+1} Y \map \tau_{\geq q} Y$.

\begin{proof}
Consider the filtration
\[
\cdots \rarr T_{\geq q+1} Y \rarr T_{\geq q} Y \rarr T_{\geq q-1 }
Y \rarr \cdots
\]
of $Y$. We have a distinguished triangle
\[
T_{\geq q+1} Y \rarr T_{\geq q} Y \rarr \Sigma^{q} \hh_{ q}(Y)
\rarr \Sigma T_{\geq q+1} Y
\]
in $\dd$ by Lemma \ref{lem:tower-fibers}. If we apply the functor $[
X , - ]_{p+q }$ and set $D_{p,q}^2 = [ X , T_{\geq q} Y]_{p +q }$
and $E_{p,q}^2 = [ X , \Sigma^{q} \hh_q (Y) ]_{p+q } $, we get an
exact couple
\[ \xymatrix{ D^2 \ar[rr]^{ (1,-1)} && D^2
\ar[dl]^{(0,0)} \\
& E^2 \ar[ul]^{(-2,1)}& } \] with the bidegrees of the maps
indicated.
 This gives a spectral sequence where
$d^r$ has bidegree $(-r , r -1)$. This follows from the definition
of the differentials given after \cite[0.6]{B}.

Now we consider conditional convergence. Recall from \cite[5.10]{B}
that we need to show that the limit $\displaystyle\lim_{p \rarr
-\infty} D_{p , n-p }^2$ and the derived limit $\displaystyle\lim_{p
\rarr -\infty}^1 D_{p , n-p }^2$ are both zero, while the map
$\displaystyle \colim_{p \rarr \infty} D_{p , n-p }^2 \rarr [X ,
Y]^n$ is an isomorphism.

For the limit and the derived limit, Lemma \ref{lim1} gives us a
short exact sequence
\[
\lim^{1}_{p \rarr -\infty }[ X , T_{\geq n-p }Y]_{n+1} \rarr [ X ,
\holim_{p \rarr -\infty } T_{\geq n-p }Y ]_n \rarr \lim_{p \rarr
-\infty } [ X , T_{\geq n-p }Y ]_n.
\]
The middle group is zero by our assumption, so the first and last
groups are also zero.
%The claim follows from Lemma
%\ref{lem:rigid-truncate}.

For the colimit, we claim that the map $\displaystyle\colim_{p \rarr
 \infty} D_{p , n-p }^2 \rarr [ X , Y ]_n$ is an isomorphism for
all $n$ if and only if $\displaystyle\colim_{q \rarr \infty} [ X ,
T_{ \leq n-q } Y]_n$ is zero for all $n$. This follows from the
distinguished triangle $T_{\geq n -p} Y \rarr Y \rarr T_{\leq n-p-1}
Y \map \Sigma T_{\geq n-p} Y$ and the fact that directed colimits of
abelian groups respect exact sequences. Under our assumption, $X$ is
weakly equivalent to $\tau_{\geq m} X$ for some $m$, so $[X, T_{\leq
n-q} Y]_n$ is zero whenever $q > -m$. Thus $\displaystyle\colim_{q
\rarr \infty} [ X , T_{ \leq n-q } Y ]_n$ is zero.
\end{proof}

We now specialize to the homotopy category $\ee$ of the
$\hh_*$-model structure on $\p \mm$. We first give a lemma which
shows that one of the conditions in Theorem \ref{specseq} is always
satisfied.

\begin{lem} \label{lem:holim}
For any $Y$ in $\p \mm$, $\holim_{q \rarr \infty} T_{\geq q } Y$ is
contractible in the $\hh_*$-homotopy category $\ee$.
\end{lem}

\begin{proof}
Each map $T_{\geq q+1} Y \map T_{\geq q} Y$ is a fibration, so the
homotopy limit is the same as the ordinary limit $\lim_{q \map
\infty} T_{\geq q} Y$. If $I$ is the indexing category for $Y$,
then one model for this limit is the pro-object $Z$ indexed by $I
\times \nn$ such that $Z_{s,q} = T_{\geq q} Y_s$ \cite[4.1]{lim}.

The map $Z_{s,q} \map *$ is an $n$-equivalence whenever $q \geq n$.
This shows that $Z \map *$ is an essentially levelwise
$n$-equivalence for all $n$, so it is an $\hh_*$-weak equivalence.
\end{proof}

\begin{thm} \label{thm:pro-specseq} Let $\mm $ be a t-model category.
Let $X$ and $Y$ be objects in $\p \mm$. There is a spectral sequence
with \[ E_{p,q}^2 = H^{-p} ( X ; \hh_q ( Y) )
\]
The spectral sequence converges conditionally to $\ee (X , \Sigma^{-p
-q } Y )$ if:
\begin{enumerate}
\item
$X$ is uniformly essentially levelwise bounded below (i.e., each
$X_s$ belongs to $\mm_{\geq n}$ for some fixed $n$), or
\item
if $Y$ is a constant pro-object and $X$ is essentially levelwise
bounded above (i.e., each $X_s$ belongs to $\mm_{\geq n}$ for some
$n$ depending on $s$).
\end{enumerate}
\end{thm}

Recall that the object $\hh_q (Y)$ by definition belongs to the
heart $\hh(\ee)$ of the $\hh_*$-homotopy category. However, when the
conditions of Proposition \ref{pro:rigid-heart} are satisfied, we
can also view $\hh_q(Y)$ as the object of $\p \hh(\dd)$ obtained by
applying $\hh_q$ to $Y$ levelwise. When $Y$ is a constant pro-object
(i.e., belongs to $\mm$, not $\p \mm$), then $\hh_q(Y)$ belongs to
$\hh(\dd)$.

\begin{proof}
The spectral sequence and conditional convergence under the first
hypothesis follows from Theorem \ref{specseq} and Lemma
\ref{lem:holim}. Observe that an object $X$ in $\p \mm$ is
bounded below (in the sense that it belongs to $\ee_{\geq n}$ for
some $n$ if and only if $X$ is uniformly essentially levelwise
bounded above; this follows from Lemma \ref{lem:pro-leq-n}.

As in the last paragraph of the proof of Theorem \ref{specseq}, it
remains to show that under the second hypothesis, $\colim_{q \rarr
\infty} \ee ( X , \Sigma^{-n} T_{\leq n-q} Y)$ vanishes for all $n$.
Because $Y$ is a constant pro-object, Corollary \ref{cor:XcY}
implies that the colimit is isomorphic to
\[
\colim_{q \rarr \infty} \colim_s \dd ( X_s , \Sigma^{-n} T_{\leq n
-q } Y).
\]
Now exchange the colimits. By hypothesis, $X_s$ belongs to
$\mm_{\geq m}$ for some $m$. Then $\dd(X_s, \Sigma^{-n} T_{\leq n-q} Y)
= 0$ for $q > -m$, so $\colim_{q \rarr \infty} \dd ( X_s ,
\Sigma^{-n} T_{\leq n -q } Y)$ vanishes for each $s$.
\end{proof}

\section{Tensor structures on pro-categories}
\label{sec:protensor}
%This section contains some results on
%tensor structures on pro-categories.
In this section we give some basic results about tensor structures
on pro-cate\-gories.

 Let $\cc$ be a tensor category. There is a \dfn{
levelwise tensor structure} on $\p \cc$ given by letting $ \{X_a \}
\otimes \{ Y_b \}$ be the pro-object $ \{X_a \otimes Y_b \}$. The
unit object of $\p \cc $ is the constant pro-object with value the
unit object in $\cc$.
 We only consider
tensor structures on $\p \cc$ that are levelwise tensor structures
inherited from a tensor structure on $\cc$.

If $\cc$ is a cocomplete category, then $\p \cc$ is a cocomplete
category \cite[11.1]{pro}. We recall the description of arbitrary
direct sums and of coequalizers in $\p \cc$.

Let $A$ be an indexing set and let $X^{\alpha} \in \p \cc$ for
$\alpha \in A$ be a set of pro-objects in $\cc$.
 Let $I_{\alpha}$ be the cofiltered
indexing category of the pro-object $X^{\alpha }$. The coproduct
$\textstyle\coprod\nolimits_{\alpha \in A} X_{\alpha}$ in $\p \mm$
is the pro-object
\[ \{ \textstyle\coprod\nolimits_{\alpha \in A}
 X_{i_{\alpha}}^{\alpha} \} \] indexed on the cofiltered category
 $ \prod_{\alpha} I_{\alpha} $.

Up to isomorphism we can assume that a coequalizer diagram is given
by levelwise maps $\{ X_a \} \rightrightarrows \{Y_a \}$.
 The coequalizer
is the pro-object $\{ \text{coeq} ( X_a \rightrightarrows Y_a ) \}$
obtained by forming the coequalizer levelwise in $\cc$.

We now consider how direct sums and tensor products interact. Let
$Y$ be a pro-object indexed on $J$. We have that
$\textstyle\coprod\nolimits_{\alpha \in A} (X^{\alpha} \otimes Y) $
is the pro-object \[ \{ \textstyle\coprod\nolimits_{\alpha} (
 X_{i_{\alpha}}^{\alpha} \otimes
Y_{j_{\alpha}} ) \} \] indexed on $
%\prod_{\alpha} ( i_{\alpha} \times j_{\alpha} ) \in
\prod_{\alpha} ( I_{\alpha} \times J ) $. On the other hand we have
that $(\textstyle\coprod\nolimits_{\alpha} X^{\alpha}) \otimes Y $
is the pro-object \[ \{ ( \textstyle\coprod\nolimits_{\alpha}
 X_{i_{\alpha}}^{\alpha} ) \otimes
Y_{j} \} \] indexed on $
% ( \prod_{\alpha} i_{\alpha} ) \times j \in
(\prod_{\alpha} I_{\alpha}) \times J $. There is a canonical map
from $\coprod_{\alpha} (X^{\alpha} \otimes Y) $ to $(
\coprod_{\alpha} X^{\alpha}) \otimes Y $.

\begin{lem} \label{smashcolimit} Let $\cc$ be a cocomplete
tensor category. If the tensor product in $\cc$ commutes with finite
direct sums (coequalizers), then the tensor product in pro-$\cc$
also commutes with finite direct sums (coequalizers).
\end{lem}
\begin{proof} This follows by cofinality arguments. \end{proof}

The tensor product on $\p \cc$ might not commute with arbitrary
direct sums even if $\cc$ is a closed tensor category. In
particular, the tensor structure on $\p \cc$ is typically not
closed.

\begin{exmp} \label{smashsum}
Let $\cc$ be a tensor category with arbitrary direct sums. Assume
that the tensor product on $\cc$ commutes with arbitrary direct
sums. Then the tensor product with a constant pro-object in $\p \cc$
respects arbitrary direct sums. In general, however, tensor product
with a pro-object does not commute with arbitrary direct sums.
 Let $X$ be a
 pro-object indexed on natural numbers.
 The canonical map $ (\coprod_{0}^{\infty} c (I)
 ) \otimes X \rarr \coprod_{0}^{\infty} (c (I)
 \otimes X ) $ is the map
\[ \left\{ \textstyle\coprod\nolimits_{i=0}^{\infty}
X_{n_i} \right\}_{ \{ n_i \} \in \mathbb{N}^{\mathbb{N}} } \rarr
\left\{ \textstyle\coprod\nolimits_{i=0}^{\infty} X_{n} \right\}_{ n
\in \mathbb{N}} .\] Assume that this map is an isomorphism, then
there is an integer $n$ and integers $n_i \geq 0$ so that $
\textstyle\coprod\nolimits_{i=0}^{\infty} X_{n_i + n+ i} \rarr
\textstyle\coprod\nolimits_{i=0}^{\infty} X_{i} $ factors through $
\textstyle\coprod\nolimits_{i=0}^{\infty} X_{n}$. Hence the map is
typically not a pro-isomorphism.

%For example, if $\cc$ has a zero object and
%if the map is an isomorphism, then there is an $n$ so that $
%so that $\textstyle\coprod\nolimits_{i=0}^{\infty} X_{i} $ is
%a retract of $\textstyle\coprod\nolimits_{i=0}^{\infty}
%X_{n}$.
%if for each $n$ we have that $X_n$ is not a retract of $X_m$
%in $\cc$ for any $m < n$.
In particular, the tensor product on the category of pro abelian
groups does not respect arbitrary sums.
\end{exmp}

 In a closed tensor category $\cc$ the tensor product with any
 object in $\cc$
  respects epic maps. The same is true
for a tensor product on $\p \cc$, even though the tensor product is
not closed.

\begin{lem} Let $\cc$ be a closed tensor category.
Then the tensor product on $\p \cc$ respects epic maps.
\end{lem} \begin{proof}
 A pro-map $f \col X \rarr Y$ is
epic if and only if for any $b$ in the indexing category $B$ of $Y$
and any two maps $Y_b \rightrightarrows Z$, for some $Z \in \cc$,
which equalize $f$ composed with the projection to $Y_b$, there is a
map $Y_{b'} \rarr Y_b$ which also equalizes the two maps to $Z$.

 Let $f \col X \rarr Y$ be an epic map. Assume that $ f \otimes
1_W $ equalizes the two maps $ Y \otimes W \rightrightarrows
 Z$ where $Z$ is an
object in $\cc$. Given two maps $ Y_b \otimes W_k \rightrightarrows
Z$. Assume that $f \otimes 1_W$ composed with the projection to $Y_b
\otimes W_k$ equalize the two maps. Then there is an $a$ and a $k'$
so that $ X_a \otimes W_{k'} \rarr Y_b \otimes W_k$ equalizes these
two maps. We have that $ X_a \rarr Y_b $ equalizes the two adjoint
maps $ Y_b
 \rightrightarrows F (W_{k'} , Z)$,
where $F$ denote the inner hom functor.
 Hence by the assumption that $f$ is epic
there is a map $Y_{b'} \rarr Y_b$ which also equalize the two maps.
Now use the adjunction one more time to get the conclusion that
$Y_{b'} \otimes W_{k'} \rarr Y_b \otimes W_{k}$
  equalize the two maps $Y_b \otimes W_{k} \rightrightarrows Z$.
\end{proof}

%\comment{We might remove the next paragraph and the following
%Lemma.}
 One might consider
monoids in pro-$\cc$. This is a more flexible notion
 than pro-objects
in the category of pro-($\cc$-monoids).
%This consists
%of pro-$\ww$-$R$-modules $X$ and pro maps $X \wedge X \rarr X
%$ and $c (R) \rarr X$ so that the usual unit and associativity
%(and commutativity) relations holds in pro-$\ww$-$\mm_R$.
We have that the category of monoids in pro-$\cc$ is the category of
algebras for the monad $\mathbb{T} X = \coprod_{n \geq 0} X^{\otimes
n} $. The category of commutative monoids in pro-$\cc$ is the
category of algebras for the monad $\mathbb{P} X = \coprod_{n \geq
0} X^{\otimes n} / \Sigma_n $.
\begin{lem} \label{smashepic}
Let $\cc$ be a complete and cocomplet closed tensor category. Then
the category of (commutative) monoids in pro-$\cc$ is complete and
cocomplete.
\end{lem}
\begin{proof} We can follow \cite[II.7]{ekmm}.
The proof of Proposition II.7.2 in \cite{ekmm} only uses that the
tensor product commutes with finite colimits and respects
epimorphisms. This holds by Lemmas \ref{smashcolimit} and
\ref{smashepic}. Hence the result follows from \cite[II.7.4]{ekmm}.
\end{proof}

\section{Tensor model categories}
\label{sec:protensormodelstructure} We give conditions that
guarantee that a tensor product on a model category $\mm$
 induces a tensor product on the homotopy
category of $\mm$. We also give more specific conditions for a
t-model category which guarantee that the induced tensor product
respects the triangulated structure and the t-structure on its
homotopy category.
%We denote the homotopy
%category of $\p \mm$ with the strict model structure by
%\mdfn{$\pp_{strict}$}.

The \dfn{pushout product axiom} for cofibrations says that if $ f
\col X \rarr Y$ and $f' \col X' \rarr Y'$ are cofibrations, then
  the \bf pushout product map \rm
\[ ( X \otimes Y' ) \amalg_{(X \otimes X' )} ( Y \otimes X' ) \rarr
( Y \otimes Y') \] is a cofibration, and if in addition $f$ or $g$
is a weak equivalence, then the pushout product map is also a weak
equivalence.
\begin{defn} A \bf tensor model category \rm $\mm$ is a model
category with a tensor product such that
\begin{enumerate} \item $\mm$ satisfies the \bf pushout product axiom
\rm
 for
 cofibrations
\item the functors $- \otimes C$ and $C \otimes -$
 take weak
 equivalences to weak equivalences for all cofibrant objects $C$ in
 $\mm$. \end{enumerate} \end{defn}

See Hovey \cite[4.2.6]{hov} for more details on tensor model
categories. Our definition is slightly stronger than his. If $\mm$
is a tensor model category, then there is a tensor product on the
homotopy category $\dd$ of $\mm$ \cite[4.3.2]{hov}. The
homotopically correct tensor product is given by first making a
cofibrant replacement of at least one of the two objects and then
form the tensor product.

 Let $\mm$ be a pointed simplicial model category and a
 symmetric
tensor category. Let $\rho $ be the functor from simplicial sets to
$\mm$ obtained by applying the simplicial tensorial structure on
$\mm$ to the unit object in $\mm$.
\begin{defn}
We say that the tensor structure and the simplicial structure on
$\mm$ are compatible if there is a natural isomorphism between the
simplicial tensorial structure and the functor $ Id \otimes \rho $,
restricted to finite simplicial complexes.
\end{defn}

\begin{lem} \label{tensorsimplicial}
Let $\mm$ be a simplicial t-model category and a tensor category. If
the tensor product and the simplicial structure are compatible, then
the simplicial structure and the levelwise tensor structure on $\p
\mm$ with the strict or $\hh_*$-model structures are also
compatible.
\end{lem}
\begin{proof} For a finite simplicial
complexes the simplicial tensorial structure on $\p \mm$ is given by
applying the simplicial structure on $ \mm$ levelwise
\cite[Sec.~16]{pro}. Hence for a finite simplicial set $K$ we have
that $ X \otimes \rho (K)$ is naturally isomorphic to the simplicial
tensor of $X$ with $K$.
\end{proof}

 We only consider the most na\"ive compatibility of the
tensor structure and the triangulated structure on the homotopy
category of a stable model category.
%A more structured version of a
%tensor triangulated category is discuss in \cite{may}.
\begin{lem}
\label{lem:triangulatedtensor}
 Let $\mm$ be a symmetric tensor
model category with a compatible based simplicial structure. Assume
the tensor product respects pushouts. Then there is a tensor
triangulated structure on the homotopy category $\dd$ in the sense
of a nonclosed version of \cite[A.2]{hps}.
\end{lem} We make use of the following property of a tensor
triangulated category in Proposition \ref{pro:AHmultiplicative}:
There are natural isomorphisms $ (\Sigma X ) \otimes Y \rarr \Sigma
(X \otimes Y)$ and $ X \otimes \Sigma Y \rarr \Sigma ( X \otimes Y)$
so the following holds: If $ X \rarr Y \rarr Z \rarr \Sigma X$ is a
distinguished triangle, then $ X \otimes W \rarr Y \otimes W \rarr Z
\otimes W \rarr \Sigma ( X \otimes W )$ is again a distinguished
triangle, where the last map is $ Z \otimes W \rarr ( \Sigma X )
\otimes W \rarr \Sigma (X \otimes W)$, and similarly when the
triangle is tensored by $W$ from the left.
\begin{proof} The unit and associativity conditions stated in
\cite[A.2]{hps} follows from the corresponding results for the
tensor product. The results stated above follows since tensor
products respects homotopy cofibers by our assumption.
\end{proof}

%The tensor structure on $\mm$ gives a tensor structure on $\p \mm$
%defined by $ \{ X_a \} \otimes \{ Y_b \} = \{ X_a \otimes Y_b \}$.

Next we give a compatibility of the tensor structure with respect to
the t-model structure. The conditions are used Proposition
\ref{pro:AHmultiplicative} to get a multiplicative structure on the
Atiyah-Hirzebruch spectral sequence.

\begin{defn} Let $\dd$ be a triangulated category.
A t-structure $\dd_{\geq 0} $ and $\dd_{\leq 0}$ and a tensor
structure on $\dd$ are compatible if $\dd_{\geq 0}$ is closed under
the tensor product.
\end{defn} Thus for all integers $i$ and $j$ we have that if $ X \in
\dd_{\geq i} $ and $Y \in \dd_{\geq j}$, then $X \otimes Y \in
\dd_{\geq i +j}$.
%Note that the axioms for a t-structure is invariant under shift in
%$n$. That is $ \dd^{\leq n} , \dd^{\geq n} $ is a t-structure on
%$\dd$ if and only if $\dd^{\leq 0} , \dd^{\geq 0} $ is a t-structure
%on $\dd$.
\begin{rem} If the t-structure on $\dd$ is not constant,
then the unit object of any tensor structure compatible with the
t-structure on $\dd$ must be an object in $ \dd_{\geq 0}$.
\end{rem}

\begin{pro} \label{tensorstrict}
Let $\mm$ be a tensor model category. Then $\p \mm $ with the strict
model structure is also a tensor model category. In particular,
there is an induced tensor structure on its homotopy category. If in
addition, the simplicial structure is compatible with the tensor
product, then the homotopy category is a triangulated tensor
category.
\end{pro}
\begin{proof}

The pushout product axiom holds since the pushout product map can be
defined levelwise.

 Let $ f $ be a weak equivalence
in $\p \mm$. We can assume that $ f$ is a levelwise weak equivalence
$ \{ f_s \col X_s \rarr Y_s \}$.
 We can furthermore assume
that the cofibrant object is a levelwise cofibrant pro-object $\{
Z_t \}$ indexed on a directed set $T$. We get that $ \{f_s \}
\otimes \{Z_t \} $ and $ \{Z_t \} \otimes \{f_s \} $ are levelwise
weak equivalences. The last statement follows from Lemmas
\ref{tensorsimplicial} and \ref{lem:triangulatedtensor}.
\end{proof}

% The
%category $\pp_{strict}$ is not a pro-category, but the tensor
%product on $\pp_{strict}$ given by Proposition
%\ref{tensorstrict} is defined using the tensor product in $\p
%\mm$ for appropriate cofibrant representatives of the objects
%in $\pp_{strict}$. In particular, the tensor product on
%$\pp_{strict}$ does not commute with direct sums in general.

We do not get an induces tensor structure on the homotopy category
$\ee$ of $\p\mm$ with the $\hh_*$-model structure when $\mm$ is a
t-model category and a tensor model category. This does not even
hold when
  the t-structure and the tensor structure on $\dd$ respect each
other. But in this case
 we do get a tensor product on the full subcategory of
 $\ee$ consisting of objects that are \ess\
bounded below.

\begin{defn} \label{defn:lessthaninfty}
Let $\mm$ be a t-model category. Let \mdfn{$\mm_{> -\infty}$ } be the
full subcategory of $\mm$ with objects $X$ so that $X \in \mm_{\geq
n}$ for some $n$. \end{defn}
 The category $\p \mm_{>
-\infty}$ is the strictly full subcategory of $\p \mm$ consisting of
objects that are \ess\ bounded below. It is larger than the category
$(\p \mm)_{> -\infty}$.
\begin{lem}
 Let $\mm$ be a t-model category.
 Then the category $\mm_{> -\infty}$ inherits a
t-model structure from $\mm$.
\end{lem}
The model structure $\mm_{>- \infty}$ has only finite colimits and
limits. The classes of cofibrations, weak equivalences, and
fibrations are all inherited from the full inclusion functor $\mm_{>
-\infty} \rarr \mm $.
\begin{proof}
If $f\col X \rarr Y$ is a map in $\mm_{> -\infty}$, and $ X
\stackrel{g}\rarr Z \rarr Y $ is a factorization of $f$ as an
$n$-equivalence followed by a co-$n$-equivalence in $ \mm$, then $Z$
is also in $\mm_{>- \infty}$: Assume that $X \in \mm_{\geq m}$ for
some $m$. In the homotopy category of $\mm$ we have a triangle $
\fib ( g) \rarr X \stackrel{g}{\rarr} Z$. Hence by Corollary
\ref{t2outof3} we have that $Z$ is in $\dd_{\geq \text{min} \{m ,n\}
}$ so $Z \in \mm_{>- \infty}$. A similar argument shows that we have
functorial factorizations of any map in $\mm_{>- \infty}$ as an
acyclic cofibration followed by a fibration and as a cofibration
followed by an acyclic fibration. The rest of the t-model category
axioms are inherited from $\mm$.
\end{proof}
%\begin{sch} In general if $\mm'$ is a full subcategory of
%of a model category $\mm$ so that $\mm'$ contains all object weakly
%equivalent to an object of $\mm'$, then $\mm'$ inherits a
%model structure from $\mm$. We need to assume that $\mm'$ has
%finite limits and colimts if we include this in the axioms of
%a model category. Note that the (co)limits in $\mm'$
%might not agree with the (co)limits in $\mm$. \end{sch}

\begin{pro} Let $\mm$ be a t-model category with a
tensor model structure
 so that the tensor product on $\dd$ is compatible with
the t-structure. Then the model category $\p \mm_{>- \infty} $ is a
tensor model category and the tensor product on its homotopy
category is compatible with the t-structure. If in addition, the
simplicial structure and the tensor structure on $\mm$ are
compatible, then $\text{Ho} \, ( \p \mm_{>- \infty} )$ is a tensor
triangulated category.
\end{pro}

\begin{proof}
Let $ f $ be a weak equivalence in $\p \mm$. We can assume that $ f$
is a levelwise map $ \{ f_s \col X_s \rarr Y_s \}$ indexed on a
directed set $S$ such that for all $n$ there is an $s_n$ such that
$f_s$ is $n$-connected for all $s \geq s_n$ \cite[3.2]{ffi}. We can
assume that the cofibrant object is a pro-object $\{ Z_t \}$ indexed
on a directed set $T$ so that $ Z_t \in \dd_{\geq n_t }$ and each
$Z_t$ is cofibrant. We use that tensoring with a cofibrant object
has the correct homotopy type. The indexing set $ \{ s , t \in S
\times T \, | \, \text{conn} (f_s) + \text{conn} ( X_t ) \geq n \}$
is cofinal in $S \times T$. Hence we
 have that $ \{f_s \} \otimes \{Z_t \} $ is an essentially levelwise
 $ n$-equivalence for all $n$.

The first part of the pushout-product axiom follows by considering
two levelwise cofibrations. When one of the maps is a levelwise
acyclic cofibration we use the previous paragraph and Lemma
\ref{proper} to show that the pushout-product map is also a weak
equivalence. The last statement follows from Lemmas
\ref{tensorsimplicial} and \ref{lem:triangulatedtensor}.
\end{proof}

\subsection{Multiplicativity in the Atiyah-Hirzebruch spectral
sequence}

We show that if $Y$ is a monoid in a tensor triangulated category
with a t-structure that is compatible with the tensor
 structure, then the Atiyah-Hirzebruch spectral sequence is
multiplicative.
%More generally, a pairing $ X \wedge Y \rarr Z$ gives a pairing of
%spectral sequences.
\begin{pro}
\label{pro:AHmultiplicative} Let $\dd$ be a symmetric tensor
triangulated category with a compatible t-structure. Let $Y$ be a
monoid in $\dd$. Then the spectral sequence in \ref{specseq} is
multiplicative.
\end{pro}
\begin{proof}
 For
convenience let $h_n$ denote $\tau_{\leq n} \tau_{\geq n} \cong
\Sigma^{n} \hh_n $.
 It suffices
to prove that we have unique dotted maps
\[ \xymatrix{ Y_{\geq i} \otimes Y_{\geq j } \ar@{.>}[r]^f \ar[d] &
Y_{\geq i + j } \ar[d]
\\ h_i ( Y) \otimes h_j (Y) \ar@{.>}[r]^g & h_{i +j } (Y) } \]
where $f$ is compatible with the multiplication on $Y$. Consider the
square \[ \xymatrix{ Y_{\geq i} \otimes Y_{\geq j } \ar[d]
\ar@{.>}[r]^f & Y_{\geq i + j } \ar[d] &
\\ Y \otimes Y \ar[r] & Y \ar[r] & Y_{\leq i+j -1 } . } \]
Since $Y_{\geq i} \otimes Y_{\geq j } \in \dd_{\geq i +j}$ we get
that the map from $Y_{\geq i} \otimes Y_{\geq j } $ to $Y_{\leq i +
j -1}$ vanish. Hence there is a lift to $Y_{\geq i +j}$. This lift
is unique since the difference of two lifts factors through
$\Sigma^{-1} Y_{\leq i + j -1} \in \dd_{\leq i+j -2}$. Hence there
is a unique map $f$. We now prove that there is a unique map $g$
between the cohomology. Consider the square
\[ \xymatrix{ Y_{\geq i+1} \otimes Y_{\geq j} \ar[r]
 & Y_{\geq i} \otimes Y_{\geq j } \ar[r] \ar[d] &
Y_{\geq i + j } \ar[d]
\\ & h_i( Y) \otimes Y_{\geq j} \ar@{.>}[r] & h_{i +j } (Y) . } \]
Since $Y_{\geq i+1} \otimes Y_{\geq j } \in \dd_{\geq i + j +1}$ and
$\Sigma Y_{\geq i+1} \otimes Y_{\geq j } \in \dd_{ \geq i+j +2}$,
 we get that there is a unique map
making the diagram commute. A similar argument with the
distinguished triangle $Y_{\geq j +1} \rarr Y_{\geq j} \rarr h_j (Y)
\rarr \Sigma Y_{\geq j+1 }$ tensored from the right by $h_i (Y)$
gives a unique map $ h_i ( Y) \otimes h_j (Y) \rarr h_{i +j } (Y) $
compatible with $Y_{\geq i } \otimes Y_{\geq j} \rarr Y_{\geq i +j }
$.
\end{proof}

%\begin{rem}
%Assume that $\dd$ is a closed tensor triangulated category. Let $F$
%denote the inner hom functor in $\dd$, and let $I$ denote the unit
%%topology we can define $E$-homology and $E$-cohomology.
 % The $n$-th $E$-homology of $X$ is $ \hh^{-n} (X \wedge
%E)$. The $n$-th $E$-cohomology of $X$ is $E^n ( X ) = \hh^n ( F ( X
%, E ))$. We have that the homology is a covariant functor, and the
%cohomology is a contravariant functor from $\dd$ to the heart $\hh
%(\dd)$. For the homotopy category of spectra the equivalence
%between the heart and the category of abelian groups gives a natural
%identifications of
% $E^* (X)$ with
%\[ \pi_0 ( E^* (X) ) \cong
% [ I , E^* (X) ] \cong [ X , E ] ^*\] and $E_* (X)$ with
% \[ \pi_0 ( E_* (X) ) \cong [ I , E_* (X) ] \cong [ I , X \wedge E ]_* . \]
%So there is no need to distinguishing between these two pairs of
%homology and cohomology functors.
%One can use the filtration on $E$ given by the t-structure to give
%spectral sequences in the heart $\hh (\dd ) $ of $\dd$.
%\end{rem}

The following remark says that it suffices to consider monoids in
the homotopy category of $\p \mm$ with the strict model structure to
get a multiplicative structure on the Atiyah-Hirzebruch spectral
sequence of Theorem \ref{thm:pro-specseq}.

\begin{rem}
We can apply Theorem \ref{specseq} to the homotopy category of $\p
\mm$ with the strict model structure and with the levelwise
t-structure. We get a spectral sequence that is isomorphic to the
spectral sequence obtained from $\ee$ with the levelwise
t-structure. This is seen by inspecting the definition of $D^2$,
$E^2$ and the differentials using Propositions
\ref{pro:mappingspace} and \ref{cor:homset}.

%This give the following description of the $\hh_*$-homotopy
%category. The homsets of the homotopy category of the
%$\hh_*$-model structure are the abutment (under reasonable
%conditions) of the Atiyah-Hirzebruch spectral sequence
%obtained from the levelwise t-structure applied to the
%homotopy category of the strict model structure on $\p \mm$.
\end{rem}

\end{document}